\documentclass[preprint,12pt]{elsarticle}

\usepackage{amsmath,amsthm,amssymb,amsfonts}
\usepackage{multirow}
\usepackage{varwidth}
\usepackage{booktabs}
\usepackage{fancyvrb}
\usepackage[ruled,linesnumbered]{algorithm2e}

\SetAlFnt{\small}
\SetAlCapFnt{\small}
\SetAlCapNameFnt{\small}
\SetAlCapHSkip{0pt}
\DontPrintSemicolon
\usepackage{graphicx}
\usepackage{caption}
\usepackage{subcaption}
\usepackage[hidelinks]{hyperref}
\usepackage{url}
\usepackage[T1]{fontenc}
\usepackage[left=2.5cm,right=2.5cm,top=3cm,bottom=3cm]{geometry}
\usepackage[capitalize]{cleveref}

\def\xref{\mathbf{x}_\mathrm{ref}}
\def\nref{\mathbf{n}_\mathrm{ref}}
\def\qtot{Q_\mathrm{tot}}
\def\ubulk{U_\mathrm{bulk}}

\journal{Applied Mathematics and Computation}

\begin{document}

\begin{frontmatter}

  \title{Acoustic modal analysis with heat release fluctuations using nonlinear eigensolvers\tnoteref{grants}}
  \tnotetext[grants]{Jose E.\ Roman was supported by the Spanish Agencia Estatal de Investigaci\'{o}n under grant PID2019-107379RB-I00 / AEI / 10.13039/501100011033.}

  \author{Varun Hiremath}
  \ead{varunhiremath@gmail.com}
  \address{Senior Software Engineer, Siemens Digital Industries Software, Simulation and Test Solutions, Simcenter Development, USA}
  \author{Jose E. Roman\corref{cor}}
  \cortext[cor]{Corresponding author.}
  \ead{jroman@dsic.upv.es}
  \address{D.\ Sistemes Inform\`atics i Computaci\'o, Universitat Polit\`ecnica de Val\`encia, Cam\'{\i} de Vera s/n, 46022 Val\`encia, Spain}

  \begin{abstract}
  Closed combustion devices like gas turbines and rockets are prone to thermoacoustic instabilities. Design engineers in the industry need tools to accurately identify and remove instabilities early in the design cycle. Many different approaches have been developed by the researchers over the years. In this work we focus on the Helmholtz wave equation based solver which is found to be relatively fast and accurate for most applications. This solver has been a subject of study in many previous works. The Helmholtz wave equation in frequency space reduces to a nonlinear eigenvalue problem which needs to be solved to compute the acoustic modes. Most previous implementations of this solver have relied on linearized solvers and iterative methods which as shown in this work are not very efficient and sometimes inaccurate. In this work we make use of specialized algorithms implemented in SLEPc  that are accurate and efficient for computing eigenvalues of nonlinear eigenvalue problems. We make use of the n-tau model to compute the reacting source terms in the Helmholtz equation and describe the steps involved in deriving the Helmholtz eigenvalue equation and obtaining its solution using the SLEPc library.
  \end{abstract}

  \begin{keyword}
    Helmholtz wave equation \sep nonlinear eigenvalue problem \sep Krylov methods \sep SLEPc
  \end{keyword}
\end{frontmatter}


\section{Introduction}\label{sec:intro}
Thermoacoustic instabilities pose a significant operational challenge and safety issue for gas turbines, rockets, and other closed combustion devices \cite{Liewen:2006:AIAA, Poinsot:2006:TNC}. These instabilities arise when there is a positive feedback between unsteady heat release oscillations in the flame and acoustic pressure waves in the chamber \cite{Rayleigh:1878}. These instabilities can cause noise, structural vibrations, flame extinction, and in worst cases even damage the device. Thus it becomes imperative to study, identify and mitigate these instabilities for safe operation.

Thermoacoustic instabilities are highly sensitive to geometry and operating conditions \cite{Juniper:2018:ARFM}, and finding instabilities late in the design process can lead to costly redesigns and delays. Design engineers seek tools for early identification and mitigation of instabilities. Several different approaches have been used for studying thermoacoustic instability. Listed here are some of the methods in decreasing order of accuracy: full scale acoustic computations using, e.g., direct numerical simulation (DNS) or large eddy simulation (LES) based methods \cite{Staffelbach:2009:PCI, Wolf:2012:FTC, Indlekofer:2017:Thesis, Xia:2019:Thesis}; methods based on linearized Navier-Stokes Equations \cite{Avdonin:2019:PCI,Bissuel:2018:CF}; solutions based on the Helmholtz equation \cite{Juniper:2018:PRF,Silva:2017:JEGTP}; and Low Order Models (LOMs) \cite{Evesque:2002:ASME,Schuermans:2003:ASME, Laurent:2019:CNF, Laurent:2021:JCP}. As we know, there is always a trade off between accuracy and computational cost. While the high end full scale simulations like LES are accurate, they are prohibitively expensive for most industrial applications, especially in the early stages of the design cycle. For this reason, the low end methods like Helmholtz solver and low-order network models are more widely used in the industry. Among these the Helmholtz equation based solvers are particularly attractive as they can be used as a post-processing tool after running the CFD simulation, so are computationally cheaper and allow to capture the geometric and boundary details more accurately \cite{Nicoud:2007:AIAA, Kaufmann:2002:CNF, Silva:2013:CNF, Falco:2021:Thesis}.

In this work, we present our implementation of the Helmholtz equation in STAR-CCM+ \cite{starccm+} software which makes use of the SLEPc \cite{Hernandez:2005:SSF} package to solve the eigenmodes to obtain the acoustic frequencies, growth rates, and the mode shapes. STAR-CCM+ is a multiphysics computational fluid dynamics (CFD) software that can be used to simulate a wide range of physics problems. In this software we have now implemented a thermoacoustic Helmholtz solver that can be used to study acoustic modes in a non-reacting or reacting simulation. This Helmholtz solver makes use the SLEPc package to solve the eigenvalue problem originating from the Helmholtz equation. SLEPc provides a wide range of algorithms to solve different kinds of eigenvalue problems which will be described in the later sections.

\section{Governing equations}
We consider a compressible reacting gas mixture. We follow the notation and derivations from \cite{Falco:2021:Thesis, Nicoud:2007:AIAA}. For acoustic analysis the viscous effects are generally neglected, and so the conservation equations for mass, momentum, and energy are given as follows
\begin{equation}
\label{eq:mass}
    \frac{D\rho}{Dt} = -\rho \nabla\cdot\mathbf{u},
\end{equation}
\begin{equation}
\label{eq:mom}
    \rho \frac{D \mathbf{u}}{Dt} = -\nabla p,
\end{equation}
\begin{equation}
\label{eq:energy}
    \rho \frac{De}{Dt} + p \nabla\cdot\mathbf{u} = \dot{q},
\end{equation}
where $\rho$ is the density, $\mathbf{u}$ is the velocity, $p$ is the pressure, $e$ is the internal energy, and $\dot{q}$ is the heat release rate.

Assuming ideal gas we get the equation of state
\begin{equation}
\label{eq:eos}
    p = \rho r T,
\end{equation}
and from thermodynamics we get
\begin{equation}
\label{eq:h-e}
    h = e + \frac{p}{\rho},
\end{equation}
where $r$ is the specific gas constant, $T$ is temperature, and $h$ is the specific enthalpy.

Taking the total derivatives of \cref{eq:eos,eq:h-e} and combining with the energy equation \cref{eq:energy} we get the following coupled pressure-energy equation (see \cite{Falco:2021:Thesis, Nicoud:2007:AIAA} for detailed derivation),
\begin{equation}
\label{eq:p-q}
    \frac{Dp}{Dt} + \gamma p \nabla\cdot\mathbf{u} = (\gamma - 1) \dot{q},
\end{equation}
where $\gamma$ is the heat capacity ratio.

\subsection{Linearization}
To study the thermoacoustic instability, we want to study the evolution of flow variables when subjected to small perturbations. To this end, we can linearize the conservation equations presented in the previous section by expressing all the instantaneous quantities as the sum of a mean and fluctuating component as follows
\begin{equation}
    \mathcal{Q}(\mathbf{x},t) = \bar{\mathcal{Q}}(\mathbf{x}) + \mathcal{Q}'(\mathbf{x},t),
\end{equation}
where $\mathcal{Q}$ denotes any quantity (e.g., pressure $p$), $\bar{\mathcal{Q}}$ denotes the time-averaged mean of that quantity, and $\mathcal{Q}'$ denotes the unsteady fluctuating component of that quantity.

We split all the flow variables $p$, $\mathbf{u}$, $\dot{q}$, etc., using the above notation, and substitute them in the conservation equations for mass \cref{eq:mass} and pressure-energy \cref{eq:p-q}. In addition we assume a low Mach number flow (which is reasonable within most combustion devices) which means $\bar{\mathbf{u}} \approx 0$. Under these assumptions, dropping all higher order fluctuations, we get the following linearized equations
\begin{equation}
\label{eq:lin-mass}
    \bar{\rho} \frac{\partial \mathbf{u}'}{\partial t} = -\nabla p',
\end{equation}
\begin{equation}
\label{eq:lin-pq}
    \frac{\partial p'}{\partial t} + \gamma \bar{p} \nabla \cdot \mathbf{u}' = (\gamma - 1) \dot{q}'.
\end{equation}

\subsection{Acoustic Wave Equation}
Combining the linearized equations \cref{eq:lin-mass,eq:lin-pq} and making approximations as explained in \cite{Poinsot:2006:TNC, Benoit:2005:IJONM, Silva:2013:CNF, Falco:2021:Thesis} gives us the following linearized acoustic wave equation
\begin{equation}
\label{wave-eqn}
  \frac{\partial^2 p'}{\partial t^2} - \nabla\cdot(\bar{c}^2 \nabla p') = (\gamma - 1)\frac{\partial \dot{q}'}{\partial t},
\end{equation}
where $p'$ is the acoustic pressure, $\bar{c} = \sqrt{\gamma  \bar{p}/\bar{\rho}}$ is the mean speed-of-sound, $\dot{q}'$ is the heat release fluctuation, $\gamma$ is the specific heat ratio, and $t$ is the time.

\subsection{Harmonic Functions}\label{sec:harmonics}
For the linearized acoustic wave equation \cref{wave-eqn} it is natural to represent $p'$ and $\dot{q}'$ in terms of complex harmonic fluctuations at frequency, $f = \omega/(2\pi)$, as follows
\begin{equation}
    p' = \hat{p} e^{- \iota \omega t},
\end{equation}
and 
\begin{equation}
    \dot{q}' = \hat{q} e^{- \iota \omega t},
\end{equation}
where $\hat{p}$ is the acoustic pressure amplitude, and $\hat{q}$ is the unsteady heat release fluctuation amplitude. Both $\hat{p}$ and $\hat{q}$ are Fourier amplitudes that are function of the frequency $\omega$.

Substituting these in the wave equation \cref{wave-eqn} we get the following Helmholtz equation in the frequency space
\begin{equation}
\label{modal-eqn}
    \omega^2 \hat{p} + \nabla\cdot(\bar{c}^2 \nabla \hat{p}) = \iota \omega (\gamma - 1)\hat{q},
\end{equation}
whose solution gives us the acoustic frequencies $\omega$ and the corresponding acoustic pressure wave amplitude $\hat{p}(\omega)$.

\subsection{Flame Transfer Function}
To solve the acoustic wave Helmholtz equation \cref{modal-eqn} and compute the acoustic frequencies and mode shapes, we need a closure model to represent the unsteady heat release rate $\hat{q}$ in terms of the pressure oscillations $\hat{p}$. The model used to express the unsteady heat release source in terms of the pressure oscillations is referred to as a Flame Response model or a Flame Transfer Function (FTF). Many models have been proposed in the past to model the heat release rate \cite{Schuller:2003:CNF,Magina:2016:CNF, Merk:2019:PCI, Meja:2018:CNF}. In this work, we use one of the simpler models referred to as the $n-\tau$ model based on the seminal work of Crocco \emph{et al.}~\cite{Crocco:1951:JARS,Crocco:1952:JARS, Crocco:1956:JFM}. This model has been successfully used in many previous studies \cite{Nicoud:2007:AIAA, Truffin:2005:CNF, Kaufmann:2002:CNF}.

\subsection{$n-\tau$ model}\label{sec:n-tau}
The $n-\tau$ model stipulates that the unsteady
heat release rate at any point in the flame is proportional to a time-lagged acoustic velocity (which is related to the acoustic pressure by \cref{eq:lin-mass}) originating from a reference upstream location.

Mathematically this model can be expressed as follows \cite{Nicoud:2007:AIAA}:
\begin{equation}
\label{n-tau-model}
\frac{\dot{q}'(\mathbf{x},t)}{\qtot} = n(\mathbf{x}) \frac{\mathbf{u}'(\xref, t - \tau(\mathbf{x}))\cdot\nref}{\ubulk},
\end{equation}
where $n(\mathbf{x})$ and $\tau(\mathbf{x})$ are fields of interaction index and time lag, respectively. $\xref$ and $\nref$ are reference location and reference direction of acoustic perturbation, respectively. $\qtot$ and $\ubulk$ are reference scaling quantities used to non-dimensionalize the relationship such that $n(\mathbf{x})$ is dimensionless.

The interaction index $n(\mathbf{x})$ relates the amplitude of heat release perturbation to acoustic velocity, and $\tau(\mathbf{x})$ estimates the time required for the acoustic perturbations to travel to the flame location. The time lag also indirectly determines the phase between the acoustic velocity and heat release, and is thus important in determining if the coupling between pressure and heat release is in-sync (unstable) or out-of-sync (stable).

The unsteady heat release source can be expressed using the $n-\tau$ model \cref{n-tau-model} as follows:
\begin{equation}
\label{qdot-n-tau}
    \frac{\partial \dot{q}'(\mathbf{x},t)}{\partial t} = n(\mathbf{x}) \left(\frac{\qtot}{\ubulk}\right) \frac{\partial}{\partial t}\left( \mathbf{u}'(\xref, t - \tau(\mathbf{x}))\cdot\nref \right).
\end{equation}

Under the low Mach number assumption, using the linearized momentum equation \cref{eq:lin-mass} we can replace $\mathbf{u}'$ with $p'$ in the heat release rate expression \cref{qdot-n-tau} which yields
\begin{equation}
    \frac{\partial \dot{q}'(\mathbf{x},t)}{\partial t} = - n(\mathbf{x}) \left(\frac{\qtot}{\rho(\xref)\ubulk}\right)  \nabla p'(\xref,t-\tau(\mathbf{x}))\cdot\nref.
\end{equation}

Substituting the harmonic relations from \cref{sec:harmonics}, we get
\begin{equation}
\label{q-n-tau-model}
\hat{q} = n(\mathbf{x}) \left(\frac{\qtot}{\iota \omega \rho(\xref) \ubulk}\right) e^{\iota \omega \tau(\mathbf{x})} \nabla \hat{p}(\xref)\cdot\nref,
\end{equation}
which is the general expression for the $n-\tau$ model and allows to close the Helmholtz equation \cref{modal-eqn} by expressing the unsteady heat release source in terms of the pressure gradient. The right hand source of \cref{modal-eqn} can thus be expressed as follows
\begin{equation}
\label{modal-rhs}
    \iota \omega (\gamma - 1) \hat{q} = n(\mathbf{x}) (\gamma - 1) \left(\frac{\qtot}{\rho(\xref) \ubulk}\right) e^{\iota \omega \tau(\mathbf{x})} \nabla \hat{p}(\xref)\cdot\nref.
\end{equation}

To use the $n-\tau$ model as given by \cref{q-n-tau-model}, we need to specify the interaction index field $n(\mathbf{x})$, the time lag field $\tau(\mathbf{x})$ and all the other reference quantities $\xref$, $\nref$, $\rho(\xref)$, $\qtot$ and $\ubulk$.

\subsubsection{Constant Time Lag}\label{sec:const-tau}
In many previous studies \cite{Nicoud:2007:AIAA, Kaufmann:2002:CNF} the time lag field $\tau(\mathbf{x})$ is assumed to be a constant in the $n-\tau$ model. This assumption allows us to significantly simplify the computation using the $n-\tau$ model. When the time lag field is constant $\tau(\mathbf{x}) = \tau$, we can express the source term in \cref{modal-rhs} as follows
\begin{equation}
\label{simplified-n-tau}
    \iota \omega (\gamma - 1) \hat{q} = S(\mathbf{x}) e^{\iota \omega \tau},
\end{equation}
where $S$ accounts for all the terms that are independent of the frequency $\omega$
\begin{equation}
\label{eq:Sx}
    S(\mathbf{x}) = n(x) (\gamma - 1) \left(\frac{\qtot}{\rho(\xref) \ubulk}\right) \nabla \hat{p}(\xref)\cdot\nref.
\end{equation}
Later in \cref{sec:nonlinear} we explain how this simplification is exploited to solve the nonlinear acoustic eigenvalue problem efficiently using SLEPc.

\subsubsection{Simplified Interaction Index}\label{sec:simp-nx}
The interaction index $n(\mathbf{x})$ needs to be specified to use the $n-\tau$ model and compute the source term $S(\mathbf{x})$ as given by \cref{eq:Sx}. The interaction index can be computed and specified in different ways. It can be computed experimentally or computationally, e.g., through an LES simulation, but this is typically very expensive and time consuming. In \cite{Nicoud:2007:AIAA}, a simpler piece-wise constant approximation of $n(\mathbf{x})$ is used (as given by Eq.~(58) in \cite{Nicoud:2007:AIAA}) which assumes $n(\mathbf{x}) = \eta$, a constant value in the flame, and $n(\mathbf{x}) = 0$ everywhere outside.

Based on this simpler piece-wise constant approximation, we use our own simplified form of interaction index, where we assume $n(\mathbf{x})$ to be proportional to a specified constant $\eta$ and the local heat release rate field as computed by the CFD simulation. So the interaction index, $n(\mathbf{x})$, is expressed as follows:
\begin{equation}
\label{eq:simp_nx}
    n(\mathbf{x}) = \eta \bar{q}_n(\mathbf{x}),
\end{equation}
where $\bar{q}_n(\mathbf{x})$ is the normalized heat release rate (as obtained from the converged CFD solution) which goes from 0 to 1, with $\bar{q}_n(\mathbf{x}) = 0$ everywhere outside the flame, and $\bar{q}_n(\mathbf{x}) > 0$ only within the flame with $\bar{q}_n(\mathbf{x}) = 1$ at the location of peak heat release rate in the flame.

The interaction index as given by \cref{eq:simp_nx} is used in this study and when presenting results we will specify the value of the interaction index coefficient $\eta$. All the remaining quantities needed to compute the source $S(\mathbf{x})$ in \cref{eq:Sx} are computed internally in STAR-CCM+ using the converged CFD solution.

\subsection{Boundary Conditions}\label{sec:bc}
The acoustic modal equation \cref{modal-eqn} requires appropriate boundary conditions. In STAR-CCM+, we have implemented the following four boundary conditions (based on the definitions from \cite{Nicoud:2007:AIAA}):
\begin{enumerate}
    \item Perfectly Reflecting (Hard Wall):
    This condition is appropriate for hard wall boundaries or inlets where the acoustic velocity fluctuation is zero. Mathematically this is imposed as
    \begin{equation}
        \nabla \hat{p}.\bf{n}_{b} = 0,
    \end{equation}
    where $\bf{n}_b$ is the unit normal vector at the boundary. 
    \item Zero Acoustic Pressure: This condition is appropriate for open boundaries with imposed pressure where the acoustic pressure fluctuation is zero giving
    \begin{equation}
        \hat{p} = 0.
    \end{equation}
    \item Constant Impedance: This condition is used to model boundary material which may partially absorb/reflect the incident wave. This property is characterized by a complex valued impedance, $Z$.
    Mathematically, this condition is imposed as follows:
    \begin{equation}
        \bar{c} Z \nabla \hat{p}\cdot\bf{n}_{b} - \iota \omega \hat{p} = 0.
    \end{equation}
    \item General Impedance: Typically the boundary impedance value can be a function of the incident wave frequency, i.e., $Z \equiv Z(\omega)$, in which case we provide a way of specifying the general impedance value $Z(\omega)$ using a simplified quadratic expression \cite{Nicoud:2007:AIAA} as follows:
    \begin{equation}
    \label{quad-imped}
        \frac{1}{Z} = \frac{1}{Z_0} + \omega Z_1 + \frac{Z_2}{\omega}.
    \end{equation}
\end{enumerate}

\subsection{Acoustic Eigenvalue Problem}\label{sec:eigenproblem}
The acoustic wave Helmholtz equation given by \cref{modal-eqn} can be solved using a Finite Element (FE) or Finite Volume (FV) based framework. In the FE approach, the quantities are computed and stored at the nodes, and the gradients are computed from shape functions \cite{Nicoud:2007:AIAA, Falco:2021:Thesis}, whereas in the FV approach, the quantities are computed and stored at the cell centers and the gradients are computed using the cell and its neighboring cell values \cite{Indlekofer:2017:Thesis}.

\Cref{modal-eqn} is solved in STAR-CCM+ using the FV formulation. For any given case, the base flow solution is first obtained on a meshed domain, following which the acoustic modal equation \cref{modal-eqn} is solved on the same mesh using the FV approach.

The first step to discretizing \cref{modal-eqn} is to integrate the equation on a cell which gives
\begin{equation}
    \int_V \omega^2 \hat{p} dV + \int_V \nabla\cdot(\bar{c}^2 \nabla \hat{p}) dV = \int_V \iota \omega (\gamma - 1) \hat{q} dV,
\end{equation}
and using the divergence theorem for the second term we get the following equation
\begin{equation}
\label{eq:modal-integral}
    \int_V \omega^2 \hat{p} dV + \int_S (\bar{c}^2 \nabla \hat{p})\cdot\mathbf{n} dS = \int_V \iota \omega (\gamma - 1) \hat{q} dV,
\end{equation}
where the first and third terms are volume integrals over the cell volume $V$, and the second term is a surface integral over the faces of that cell. The boundary conditions are applied through the second term when integrated over the boundary faces.

The above equation \cref{eq:modal-integral} is discretized and assembled on a CFD mesh resulting in an  eigenvalue problem of the following form:
\begin{equation}
\label{eq:acoustic-evp}
    \left[A + \omega B + \omega^2 C \right] \hat{P} = D(\omega) \hat{P},
\end{equation}
where $A$, $B$, $C$ and $D$ are sparse matrices of size $N \times N$ with $N$ being the number of cells in the mesh, $\hat{P}$ is a size $N$ eigenvector composed of acoustic pressure $\hat{p}$ values from all the cells, and $\omega$ the acoustic frequency is the eigenvalue obtained from the above equation. We solve this equation using the SLEPc package.

Additional details on the matrices appearing in \cref{eq:acoustic-evp} are:
\begin{itemize}
    \item $A$: This matrix accounts for the wave propagation term $\bar{c}^2 \nabla \hat{p}$ appearing in the second surface integral term in \cref{eq:modal-integral}. This matrix also gets contributions from (non-impedance) boundary conditions.
    \item $B$: This matrix accounts for the impedance boundary condition contributions originating from the wave propagation term $\bar{c}^2 \nabla \hat{p}$ appearing in the surface integral term in \cref{eq:modal-integral} on the boundaries. If none of the boundaries use the impedance boundary condition then this matrix is zero.
    \item $C$: This matrix comes from the first volume integral term in \cref{eq:modal-integral} involving the cell center $\hat{p}$ values. This matrix is identity except when using the general impedance boundary condition as given by \cref{quad-imped}.
    \item $D$: This matrix accounts for the heat release source on the right-hand-side of \cref{eq:modal-integral} and is obtained from the Flame Transfer Function $n-\tau$ model as described in \cref{sec:n-tau}. In non-reacting cases $D$ matrix is zero, and when the time lag in the $n-\tau$ model is constant then as described in \cref{sec:const-tau} this matrix can be expressed as
    \begin{equation}
    \label{eq:D-S}
        D(\omega) = e^{\iota \omega \tau} S.
    \end{equation}
\end{itemize}

\subsection{Secondary Gradient}\label{sec:second-grad}
To discretize the wave propagation term $\bar{c}^2 \nabla \hat{p}$ (appearing in the surface integral term in \cref{eq:modal-integral}) we need to compute the acoustic pressure gradient, $\nabla \hat{p}$, at the cell faces. Expressing this gradient $\nabla \hat{p}$ at the face using the adjacent cell center values $\hat{p}$ becomes challenging on an unstructured mesh.

\begin{figure}
    \centering
    \includegraphics[width=400pt]{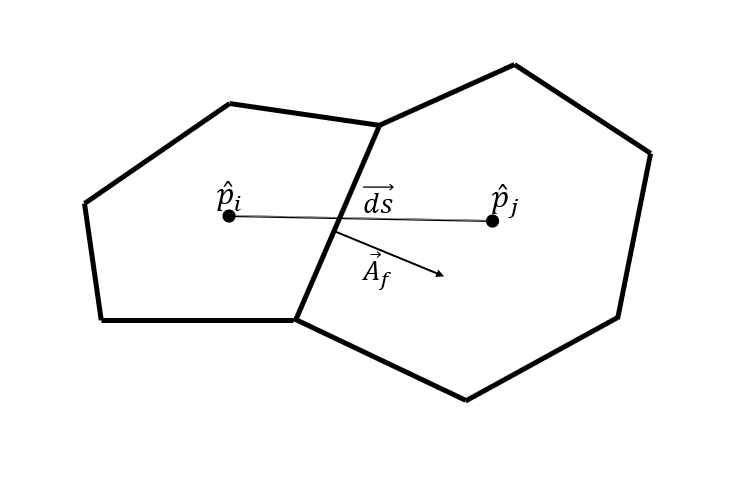}
    \caption{Schematic showing two adjacent cells in an unstructured mesh. The cell centers are separated by a distance $\vec{ds}$ and the face separating the cells has an area vector represented as $\vec{A}_f$.}
    \label{fig:gradients}
\end{figure}

If we consider two adjacent cells $i$ and $j$ from an unstructured mesh as shown in the sketch in \cref{fig:gradients}, then the acoustic pressure gradient $\nabla \hat{p}$ at the face can be expressed to a first order approximation using just the cell center values as
\begin{equation}
\label{eq:grad-1}
    \nabla \hat{p}_f = (\hat{p}_j - \hat{p}_i)\frac{\vec{A}_f}{\vec{ds}\cdot\vec{A}_f}.
\end{equation}
This however is not very accurate as in an unstructured mesh the cell faces can be highly skewed (i.e., not perpendicular to the line joining the cell centers) as seen in \cref{fig:gradients}. To correct for this skewness, a second order correction to the face gradient is applied using the gradients computed at the cell centers. There are different ways of correcting this gradient \cite{Tsui:NHT:2006, Mathur:NHT:1997, Moukalled:2016:Springer} but is typically expressed as follows
\begin{equation}
\label{eq:grad-2}
    \nabla \hat{p}_f = (\hat{p}_j - \hat{p}_i)\frac{\vec{A}_f}{\vec{ds}\cdot\vec{A}_f} + \overline{\nabla \hat{p}} - (\overline{\nabla \hat{p}}\cdot\vec{ds})\frac{\vec{A}_f}{\vec{ds}\cdot\vec{A}_f},
\end{equation}
where $\overline{\nabla \hat{p}}$ denotes the mean gradient at the face computed as follows
\begin{equation}
    \overline{\nabla \hat{p}} = \frac{\nabla \hat{p}_i + \nabla \hat{p}_j}{2},
\end{equation}
and the cell center gradients $\nabla \hat{p}_i$ are typically computed using Green-Gauss or Least-Squares method.

If we represent the acoustic pressure gradient $\nabla \hat{p}$ at the faces using the first order approximation as given by \cref{eq:grad-1}, then the sparse matrix $A$ appearing in the eigenvalue problem \cref{eq:acoustic-evp} can be explicitly built. However, if the gradient $\nabla \hat{p}$ is represented using the second order expression as given by \cref{eq:grad-2} then it is not possible to explicitly form the matrix $A$ as the expression involves the mean gradient term $\overline{\nabla \hat{p}}$ which cannot be computed without knowing the $\hat{p}$ values at the cell centers which is the unknown (eigenvector) we are trying to solve in \cref{eq:acoustic-evp}.

In the SLEPc package to solve an eigenvalue problem of the form \cref{eq:acoustic-evp} we can provide all the matrices explicitly if known. If however, a matrix cannot be explicitly built (or in some cases if storage memory is a concern), then we also have the option of providing a method which can compute matrix-vector product for any given vector. The SLEPc package refers to this as a {\it matrix-free} option provided through a special {\it shell matrix} that defines the matrix-vector product. So when representing $\nabla \hat{p}$ using the second order expression as given by \cref{eq:grad-2}, we do not explicitly build the matrix $A$. Instead we use SLEPc's {\it shell matrix} option and provide a method to compute $A \times v$ for any given vector $v$. The matrix-vector product $A\times v$ can be easily computed using our FV approach as for the given vector $v$ we can compute the gradients at the cell centers and thus obtain the mean gradient $\overline{\nabla v}$ at the face for any given vector $v$.

This {\it matrix-free} approach using the {\it shell matrix} allows us to accurately represent the face gradients, however this also poses a problem for computing the eigenvalues using SLEPc as the {\it shell matrix} cannot be inverted or factored, which in turn significantly slows down the computation. So to speed up the computation, we solve a preconditioned system where we use the matrix $A$ explicitly formed using the first order gradient given by \cref{eq:grad-1} as a preconditioner. This idea is explained later in more detail in \cref{sec:precond}.

\subsection{Eigenvalue Problem Types}
Depending on the simulation type (reacting/non-reacting) and the boundary conditions, we end up with one of the following three eigenvalue problem types:
\begin{enumerate}
    \item Non-reacting, with no impedance boundary condition. In this case we have $B = 0$ and $D = 0$ and we end up with a Linear Eigenvalue Problem of the form: 
    \begin{equation}
    \label{eq:eiglin}
        [A + \omega^2 C] \hat{P} = 0.
    \end{equation}
    \item Non-reacting, with impedance boundary condition. In this case we have $D = 0$ and we end up with a Quadratic Eigenvalue Problem of the form: 
    \begin{equation}
    \label{eq:eigquad}
        [A + \omega B + \omega^2 C] \hat{P} = 0.
    \end{equation}
    \item In the most general reacting case with impedance boundary condition we get the full Nonlinear Eigenvalue Problem of the form:
    \begin{equation}
    \label{eq:eignonD}
        [A + \omega B + \omega^2 C] \hat{P} = D(\omega) \hat{P}.
    \end{equation}
    However when using the $n-\tau$ model if $\tau$ is constant then as given by \cref{eq:D-S} we can simplify the above equation to the below form
    \begin{equation}
    \label{eq:eignonS}
        [A + \omega B + \omega^2 C] \hat{P} = e^{\iota \omega \tau} S \hat{P},
    \end{equation}
    where $S$ is a constant matrix independent of $\omega$ which helps us solve this equation efficiently using SLEPc's Nonlinear Eigenvalue Problem (\texttt{NEP}) module as explained more in the later sections.
\end{enumerate}

\section{Linearized iterative solvers} \label{sec:lin_iter}
There exist many different algorithms and solvers that can be used to compute the eigenvalues of the linear eigenvalue problem as given by \cref{eq:eiglin}. For this reason, in many of the previous implementations of the Helmholtz solver \cite{Nicoud:2007:AIAA, Falco:2021:Thesis}, the quadratic and nonlinear eigenvalue problems are solved by linearizing these equations. For instance as described in \cite{Nicoud:2007:AIAA}, the quadratic eigenvalue problem \cref{eq:eigquad} can be linearized and solved as follows
\begin{equation}
\label{eq:quad_lin}
    \begin{bmatrix}
    0 & -I\\
    A &  B
    \end{bmatrix}
    \begin{bmatrix}
    \hat{P} \\
    \omega \hat{P}
    \end{bmatrix}
    + \omega \begin{bmatrix}
    I & 0 \\
    0 & C 
    \end{bmatrix}
    \begin{bmatrix}
    \hat{P} \\
    \omega \hat{P}
    \end{bmatrix} = 0,
\end{equation}
which converts the linear eigenvalue problem given by \cref{eq:eigquad} of size $N\times N$ into a linear problem of size $2N \times 2N$. But as we will explain in the later sections, this is not the most efficient way of solving the quadratic eigenvalue problem.

Similarly, as described in \cite{Nicoud:2007:AIAA, Falco:2021:Thesis}, the nonlinear eigenvalue problem as given by \cref{eq:eignonD} can be linearized and solved iteratively. One of the algorithms as described in \cite{Nicoud:2007:AIAA} is as follows:
\begin{enumerate}
    \item As a first step, \cref{eq:eignonD} is solved by dropping the right-hand-side source involving the nonlinear term $D(\omega)$. Let one of the obtained eigen frequencies be $\omega^i$.
    \item Then a simplified (linear or quadratic) form of \cref{eq:eignonD} is solved iteratively in which the nonlinear term $D(\omega)$ is computed explicitly at the previous known value of $\omega^i_{k-1}$ at each iteration $k = 1$ to $k = n$ until convergence. The equation that is solved iteratively is given as follows
    \begin{equation}
    \label{eq:quad_iter}
        \left[(A - D(\omega^i_{k-1})) + \omega^i_k B + (\omega^i_k)^2C \right] \hat{P} = 0,
    \end{equation}
    where the iteration starts at $k = 1$ with $\omega^i_0 = \omega^i$ obtained in step 1, and continues until $\omega^i_k$ converges to within a specified tolerance. The above quadratic iterative equation \cref{eq:quad_iter} is further linearized as described above using \cref{eq:quad_lin} and then solved using an eigenvalue solver.
\end{enumerate}
The above linearized iterative algorithm works, but as explained in later sections and also reported in some previous studies \cite{Buschmann:2020:JEGTP} is not the most efficient way of solving the nonlinear equation. The above algorithm needs to be run for $n$ iterations for each eigen frequency, which makes it computationally very expensive. Moreover, for some cases where the eigen frequencies are close, this iterative algorithm can miss certain frequencies as reported in \cite{Buschmann:2020:JEGTP}. 

Both the quadratic and nonlinear eigenvalue problems can be solved more efficiently and accurately using the algorithms available in the SLEPc package as described in the next section.

\section{Eigenvalue problems and solvers}\label{sec:solvers}

The matrices involved in the eigenvalue problems discussed in the previous sections are large, sparse, and complex non-Hermitian, and we are interested in computing a selected part of the spectrum, in particular a few smallest eigenvalues $\omega$. We address the solution of these problems by using the SLEPc library~\cite{Hernandez:2005:SSF}. SLEPc provides a collection of solvers for each of the mentioned eigenvalue problems, and it is designed in a way that the user can easily switch from one solver to another for easy comparison. In this paper, we focus on Krylov-type methods because they have shown better performance and reliability for the application at hand. In this section, we give an overview of these methods, then in \cref{sec:implem} we discuss details of how these methods are used in SLEPc.

\subsection{Linear eigenvalue problems}\label{sec:linear}

In the linear eigenvalue problem, \cref{eq:eiglin}, the eigenvalue parameter is the square of the acoustic frequency, $\omega^2$. Krylov eigensolvers~\cite{Saad:2011:NML} are based on building Krylov sequences associated with a matrix $M$ and a given initial vector, $\{v_1, Mv_1, M^2v_1, \dots, M^kv_1, \dots\}$. The subspace spanned by these vectors (the Krylov subspace) will contain increasingly better approximations of the eigenvectors associated to extreme eigenvalues (particularly those of largest magnitude). The successive powers of $M$ are not computed explicitly, instead a matrix-vector product with $M$ is done at each step of the method. The generated vectors must be mutually orthonormalized to make the method numerically stable. Apart from this, the main issue to take into account when implementing these methods is what happens when slow convergence requires a large number of steps $k$, which is addressed by restart techniques as discussed below.

In our case, we have two matrices, $A$ and $C$, so it is necessary to invert one of them to be able to apply the above technique. We could either apply the Krylov method to matrix $M=-C^{-1}A$, in which case we obtain the (largest magnitude) eigenvalues of \cref{eq:eiglin}, or alternatively apply the Krylov method to matrix $M=-A^{-1}C$, in which case we get the reciprocals of the (smallest magnitude) eigenvalues of \cref{eq:eiglin}. The latter is most interesting for our application, since we need the lower frequencies. A generalization is the shift-and-invert spectral transformation, in which given a target value $\sigma$ we apply the solver to the transformed problem
    \begin{equation}\label{eq:sinvert}
        \left[(A + \sigma C)^{-1} C + \theta I\right] \hat{P} = 0,
    \end{equation}
obtaining approximations to the largest eigenvalues $\theta$ which correspond to eigenvalues of the original problem that are closest to the target $\sigma$. In other words, for computing the lower frequencies we use shift-and-invert with $\sigma=0$. The inverse $(A + \sigma C)^{-1}$ must be handled implicitly, typically by a sparse LU factorization (computed once at the beginning) and the corresponding triangular solves (every time $M$ is needed during the eigensolution).

The procedure for expanding the Krylov subspace only requires a matrix-vector product with matrix $M$. This means that $M$ need not be built explicitly, and we can instead use a subroutine for the matrix-vector product. This is usually called the matrix-free approach, or \emph{shell} matrix in PETSc's terminology. But things get complicated when the shift-and-invert strategy is employed, because shell matrices cannot be factorized, so instead of the LU decomposition for the linear solves, one has to use an iterative linear solver with a good preconditioner. In our case, the preconditioner will be the LU factorization of a first-order approximation of $A$, as discussed in \cref{sec:precond}.

The cost of orthogonalization grows with the number of steps $k$, so when a maximum size $k$ is reached without having all wanted eigenvalues converged, a restart mechanism must be invoked. A very effective restart is the Krylov-Schur method~\cite{Stewart:2001:KAL}, which consists in compressing the available $k$-dimensional Krylov subspace to a smaller dimension, such as $k/2$ for instance, so that the compressed subspace retains the most wanted approximations and purges eigendirections for unwanted eigenvalues. The compression is done in a way that the Krylov method can be continued to expand the subspace again up to the maximum dimension. This procedure is repeated until enough wanted eigenpairs have converged to the requested accuracy.

\subsection{Polynomial eigenvalue problems}\label{sec:poly}

As opposed to \cref{eq:eiglin}, the problem in \cref{eq:eigquad} is nonlinear in the eigenvalue parameter because both $\omega$ and $\omega^2$ appear in the equation. It is a quadratic eigenvalue problem~\cite{Tisseur:2001:QEP}, a particular case of the polynomial eigenvalue problem in which the matrix polynomial $P(\omega)=A + \omega B + \omega^2 C$ has degree $d=2$.

One possible approach for solving the polynomial eigenproblem is to apply a Krylov method to a linearization, that is, a linear eigenproblem of size $d\cdot N$ whose solution is directly related to the solution of the polynomial problem. This linearization technique can be applied for matrix polynomials of arbitrary degree, even in the case of polynomials expressed in non-monomial bases~\cite{Campos:2016:PKS}. In the simpler case of a quadratic eigenproblem, we can use the first companion linearization,
\begin{equation}
\label{eq:companion}
\left(\begin{bmatrix}0 & I\\-A & -B\end{bmatrix}-\omega\begin{bmatrix}I & 0\\0 & C\end{bmatrix}\right)y=0.
\end{equation}
The $2N$ eigenvalues $\omega$ of this linear eigenproblem are the same as those of the quadratic eigenproblem, and the corresponding eigenvectors can be expressed as
\begin{equation}
\label{eq:linevec}
y=\begin{bmatrix}\hat{P}\\\omega\hat{P}\end{bmatrix},
\end{equation}
where $\hat{P}$ is the eigenvector of the quadratic eigenproblem.

Implementing the above approach naively may result in inefficient calculation and possibly large numerical error. SLEPc provides a solver with all the ingredients necessary for a robust and efficient solution of polynomial eigenproblems~\cite{Campos:2016:PKS}, and we discuss below those that are more relevant for our application.

When operating on the linearization, the Krylov method has to orthogonalize vectors of length $d\cdot N$, with the corresponding increase in computational cost. For $d=2$ this is not too critical, but still it is possible to avoid this overhead. The trick is to exploit the structure of the eigenvectors, \cref{eq:linevec}, and use a compact representation of the Krylov basis in which a basis of vectors of length $N$ are used to represent both the upper and lower parts of the Krylov vectors. This results in almost half memory requirements, together with a reduction of the computational cost if the orthogonalization procedure is adapted to this representation. The resulting method is called TOAR~\cite{Lu:2016:SAT}. The next step is to adapt the Krylov-Schur restart to the new compact representation~\cite{Campos:2016:PKS}.

In polynomial eigenproblems it is also possible to apply the shift-and-invert transformation. We can either transform the polynomial and then linearize, or alternatively do the shift-and-invert transform on the linearization. The latter approach involves the LU factorization of a $2N\times 2N$ matrix, but it can be replaced with a block LU factorization in which only a factorization of an $N\times N$ matrix is required~\cite{Campos:2016:PKS}. In particular, the matrix $P(\sigma)$ is the one that has to be factorized, for $\sigma=0$ in our case, $P(0)=A$.

In terms of numerical error, the conditioning of the linearization can be quite bad in some cases, particularly when the norms of $A$, $B$ and $C$ vary wildly. In that case, a parameter scaling~\cite{Fan:2004:NSS} can be used to obtain a well-conditioned linear problem. It consists in a transformation on the eigenvalue parameter of the polynomial eigenproblem, $\omega=\rho\theta$, with $\rho:=\sqrt{{\|A\|_2}/{\|C\|_2}}$, resulting in an equivalent eigenproblem with matrix polynomial $\tilde P(\theta):=P(\rho\theta)$, with coefficient matrices $A$, $\rho B$, $\rho^2 C$, which has the same eigenvectors as the original polynomial eigenproblem, and related eigenvalues $\theta=\frac{\omega}{\rho}$.

Finally, for the extraction of the eigenvector of the quadratic eigenproblem $\hat{P}$ from the eigenvector of the linearization \cref{eq:linevec} one can normalize either the upper or the lower block, depending on the magnitude of $\omega$, otherwise severe loss of significant digits may occur.

\subsection{Nonlinear eigenvalue problems}\label{sec:nonlinear}

In the general nonlinear eigenproblem, \cref{eq:eignonD}, applying the linearization technique of the previous section will not get rid of the nonlinearity present in $D(\omega)$. We assume that the entries of the parameter-dependent matrix $D(\omega)$ can contain any nonlinear function, such as exponential, square root, etc., so we must apply numerical methods that can cope with this general nonlinear case. In recent years, significant progress has been made in this direction~\cite{Guttel:2017:NEP}.

We express the nonlinear eigenproblem of \cref{eq:eignonD} as
    \begin{equation}\label{eq:eignon2}
        F(\omega)  \hat{P} = 0,
    \end{equation}
where $F(\omega)=A + \omega B + \omega^2 C-D(\omega)$. The methods to solve this problem rely on evaluating $F(\omega)$ for various values of $\omega$, and in some methods it is also necessary to evaluate the derivative $F'(\omega)$. Since $F(\omega)$ is not a constant matrix, one way to represent it in a computer program is to provide a subroutine that, given the parameter $\omega$, evaluates the entries of $F(\omega)$, and similarly for the derivative. An alternative, usually called the split form, is to express the matrix as a sum of constant matrices multiplied by scalar functions,
    \begin{equation}\label{eq:split}
        F(\omega) = \sum_{i=0}^{\ell-1}A_if_i(\omega).
    \end{equation}
In our case, the first three functions are $f_i(\omega)=\omega^i$, $i=0,1,2$, with $A_0=A$, $A_1=B$, $A_2=C$, and finally we have $f_3=-e^{\iota \omega \tau}$ with $A_3=S$ according to the definition of $D(\omega)$ in \cref{eq:eignonS}.

SLEPc, in its \texttt{NEP} module~\cite{Campos:2021:NEP}, provides support for both ways of representing the matrix, and contains several methods to compute the solution. In this work, we use the NLEIGS method~\cite{Guttel:2014:NCF}, a Krylov-type method that work very well, especially in cases where $F(\omega)$ has singularities. This method relies on evaluations of $F(\omega)$ only, and the derivative is not required. We next give an overview of how it works.

In summary, NLEIGS applies a Krylov iteration to a companion-type linearization of a rational interpolant of the nonlinear function. The first step is to compute a rational matrix
    \begin{equation}\label{eq:interpolant}
        R(\omega) = \sum_{i=0}^{d-1}R_iq_i(\omega),
    \end{equation}
where $q_i(\omega)$ are scalar rational functions, such that $R(\omega)$ interpolates $F(\omega)$ at a number of points $\sigma_i$ located at the boundary of a region of the complex plane in which we want to search for eigenvalues. The $q_i(\omega)$ rational functions are computed from the target points $\sigma_i$ as well as points belonging to the singularity set of $F(\omega)$, while the constant matrix coefficients $R_i$ are obtained by a certain recurrence involving the evaluation of $F(\cdot)$ at the points $\sigma_i$. The number of required terms $d$ (the degree) of the rational matrix $R(\omega)$ will depend on the function being interpolated. In principle, we can assume that the solution of the rational eigenproblem $R(\omega)\hat{P}=0$ is a good approximation of the solution of \cref{eq:eignon2}.

The next step is to build a special linearization, in the spirit of \cref{eq:companion}, but with matrices of order $d\cdot N$. Again, these matrices must not be built explicitly, especially in case $d$ is large, and instead matrix-vector products are handled implicitly following the block structure of these matrices~\cite{Campos:2021:NEP}. Then, it is possible to adapt the TOAR method to this linearization, including the Krylov-Schur restart and also the shift-and-invert transformation in which the LU factorization is carried out in a blocked fashion, where the most computationally expensive operation is to factorize the $N\times N$ matrix $R(\sigma)$~\cite{Campos:2021:NEP}.

\section{Implementation details}\label{sec:implem}

In this section we discuss a few aspects related to the implementation. In particular, we describe how to organize code in situations where the problem matrices are different from the matrices used to build the preconditioner, such as the case when one of the problem matrices is implicit (not built explicitly). During this work, we have added new functionality in SLEPc to support this use case.

\subsection{Overview of SLEPc}\label{sec:slepc}

SLEPc, the Scalable Library for Eigenvalue Problem Computations~\cite{Hernandez:2005:SSF,Roman:2022:SUM}, is a parallel library that provides solvers for different classes of eigenvalue problems, including linear (both Hermitian and non-Hermitian), polynomial and general nonlinear. The solvers can be used with either real or complex arithmetic, depending on whether the problem matrices are real or complex. The parallelism model is based on MPI message passing for distributed memory computers, but it is also possible to exploit thread parallelism on the CPU (via multi-threaded BLAS for the innermost computations) or on the GPU by enabling CUDA support during installation.

The user interface for the solvers is very flexible, and it allows configuration of many settings, either in the source code or at the command-line when the program is executed. For instance, the user can easily select how many eigenvalues must be computed, which is the part of the spectrum of interest, which is the maximum dimension of the Krylov subspace, or how much accuracy is required (by specifying a tolerance).

SLEPc relies on PETSc, the Portable Extensible Toolkit for Scientific Computation~\cite{Balay:2022:PUM}, and can be seen as an extension that complements PETSc with all the functionality necessary to solve eigenproblems. PETSc provides an object-oriented hierarchy of data classes (most notably matrices and vectors, but also other related to discretization of PDE's) and solvers for systems of linear and nonlinear equations, among other. SLEPc makes heavy use of PETSc's data structures for the implementation of eigensolvers, but it also requires linear solvers in many situations, as in the case of the shift-and-invert spectral transformation discussed in \cref{sec:linear}. The functionality offered by PETSc for linear systems of equations covers both direct and iterative methods (in combination with a preconditioner). For the former, PETSc includes just a basic sequential LU (and Cholesky), so in case one wants to use a direct linear solver in parallel runs it is necessary to configure PETSc with an external package such as MUMPS \cite{Amestoy:2000:MUMPS, Amestoy:2001:FAM, Amestoy:2006:PC}. This is the recommended approach for eigenvalue computations involving shift-and-invert, unless the problem size is so large that factorization is not viable.

Apart from the matrix-vector products (or linear solves in the case of shift-and-invert), the most time-consuming part of the Krylov eigensolvers is typically the orthogonalization of vectors. This operation is critical from the numerical point of view, since eigensolvers are usually very sensitive to poor orthogonality of the basis vectors, compared to other linear algebra problems. SLEPc implements an iterated classical Gram-Schmidt orthogonalization scheme~\cite{Hernandez:2007:PAE} that is numerically robust and has high efficiency both in terms of MPI (small number of reductions) and sequentially (high arithmetic intensity).

\subsection{Preconditioning}\label{sec:precond}

As discussed in \cref{sec:eigenproblem}, matrix $A$ is related to the wave propagation term $\bar{c}^2 \nabla \hat{p}$. However when representing the face gradients using the second order expression as explained in \cref{sec:second-grad}, this matrix is not assembled explicitly. Instead, $A$ is provided as a {\it shell matrix}. This is a serious difficulty in terms of computation, because in many of the solution schemes described in \cref{sec:solvers} it is necessary to invert $A$, which is usually done by means of an LU factorization. As already mentioned, shell matrices cannot be factorized, and we must have recourse to an iterative linear solver such as GMRES or BiCGStab, provided by PETSc. However, to be practical, these iterative solvers must be combined with a good preconditioner, otherwise the number of iterations required for convergence will be exceedingly high, and ruin the performance of the overall computation.

Recall that, when solving the linear system $Ax=b$, the preconditioned system
\begin{equation}\label{eq:precond}
M^{-1}Ax=M^{-1}b
\end{equation}
is solved instead, where $M$ is chosen in a way that the number of iterations required by the iterative solver is reduced significantly. In practice, at each iteration of the method the preconditioner $M^{-1}$ is applied to the vector resulting from multiplication by $A$. A good preconditioner is one that approximates $A$ in some way, more precisely the eigenvalues of $M^{-1}A$ should all be reasonably close to 1 and far away from 0. In our case, as explained in \cref{sec:second-grad}, we use the matrix $A$ explicitly formed using the first order gradient  as a preconditioner, and use an LU decomposition of this matrix (computed by MUMPS) as the preconditioner $M^{-1}$.

In SLEPc, the matrix to be used to build the preconditioner can be provided by the user via the interface function \texttt{STSetPreconditionerMat}. However, this approach has the limitation that the matrix is constant and therefore it may not work well whenever the shift-and-invert transformation of \cref{eq:sinvert} is used with a target value $\sigma$ different from zero, as the matrix to be inverted is $A+\sigma C$. During this work, we have implemented in SLEPc a more flexible scheme, with a new user-interface function \texttt{STSetSplitPreconditioner} to provide approximations of the two matrices, $\tilde A$ and $\tilde C$, so that the preconditioner is built from $\tilde A+\sigma \tilde C$, whatever the value of $\sigma$ is.

In polynomial eigenproblems we have a similar situation. As discussed in \cref{sec:poly}, when using shift-and-invert the polynomial eigensolver needs to solve linear systems with the coefficient matrix $P(\sigma)$, which is equal to $A$ if $\sigma=0$. But for other values of $\sigma$ the approximation of $A$ alone may not be a sufficiently good preconditioner. Thus, we can use the same interface function described above to pass approximations to all the coefficient matrices of the matrix polynomial, that is, $\tilde A$, $\tilde B$ and $\tilde C$ in our application.

In the case of the nonlinear eigenproblems, this topic is even more relevant. From \cref{sec:nonlinear}, we know that the solver needs to invert matrix $R(\sigma)$, but the rational matrix $R(\cdot)$ is built internally and hence not known to the user. However, if $R(\cdot)$ is a good rational approximation to the nonlinear function, we can expect that $R(\sigma)\approx F(\sigma)$. For the case $\sigma=0$, we can see that $F(0)=A-S$ and hence the preconditioner should not be built from $A$ alone. Analogously to the previous cases, we have added a function \texttt{NEPSetSplitPreconditioner} in which the user may pass a list of $\ell$ matrices $\tilde A_i$ to build the preconditioner with an expression similar to \cref{eq:split} irrespective of the value of the argument $\omega$. In practice, the NLEIGS solver will use these approximate matrices to build an approximate form of $R(\sigma)$ that will be factorized and used as a preconditioner.

\section{Computational results}\label{sec:results}
Here we present acoustic modal solver results obtained in STAR-CCM+ using SLEPc for five different test cases. In each case, we first solve the base flow in STAR-CCM+ using the available models on a generated CFD mesh. Once we get a converged solution, we then use our Helmholtz solver to compute the acoustic modes using the base solution on the same CFD mesh. So if there are $N$ cells in the CFD mesh, then we end up with a size $N$ eigenvalue problem, where the eigenvector $\hat{P}$ is of size $N$, and all the matrices involved in the eigenvalue equation \cref{eq:eignonS} are of size $N \times N$.

In each case, we need to set appropriate boundary conditions on every boundary as described in \cref{sec:bc}. We then discretize the Helmholtz equation \cref{modal-eqn}, apply the boundary conditions and use the $n-\tau$ model to obtain the eigenvalue equation of the form \cref{eq:eignonS}. We then solve this equation using SLEPc to get the complex valued acoustic frequencies (eigenvalues) $\omega$ and the corresponding acoustic mode shapes (eigenvectors) $\hat{P}$. 

Each complex valued acoustic frequency $\omega$ has a real and imaginary part which gives us the acoustic frequency and the growth rate, respectively. Mathematically we express
\begin{equation}
    \omega = \omega_r + \iota \omega_i,
\end{equation}
where $\omega_r$ and $\omega_i$ are the real and imaginary parts of $\omega$. Substituting this in the harmonic pressure fluctuation expression as given in \cref{sec:harmonics} we get
\begin{equation}
    p' = \hat{p}e^{-\iota(\omega_r + \iota \omega_i)t} = \left[\hat{p}e^{\omega_i t} \right] e^{-\iota \omega_r t},
\end{equation}
which means if 
\begin{itemize}
    \item $\omega_i > 0$, then the acoustic fluctuations will grow in time, i.e., an {\it unstable} mode
    \item $\omega_i < 0$, then the acoustic fluctuations will subside in time, i.e., a {\it stable} mode
    \item $\omega_i = 0$, then we have a standing wave
\end{itemize}
and in addition $\omega_r$ gives us the acoustic frequency at which the mode is excited.

So for each acoustic frequency $\omega$ that we obtain from SLEPc, we extract the real (frequency) and imaginary (growth rate) parts to identify the stable/unstable modes. The frequency is typically reported as $f_r = \omega_r/(2\pi)$ in Hz, while the growth rate is expressed either as $\omega_i$ in rad/s or as $f_i = \omega_i/(2\pi)$ in Hz. For each mode the acoustic mode shape is displayed as the normalized eigenvector $\hat{P}$ obtained from SLEPc. The normalization scales the eigenvector such that the acoustic pressure is bounded between 0 and 1.

\subsection{Helmholtz Resonator}\label{sec:resonator}
\begin{figure}
    \centering
    \includegraphics[width=200pt]{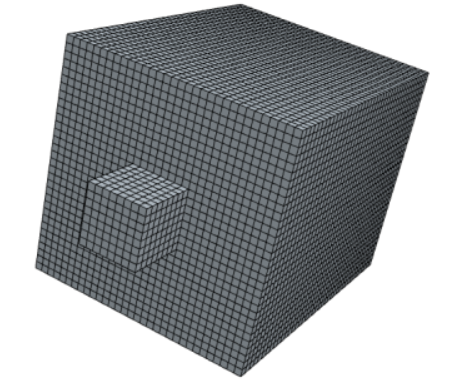}
    \caption{The CFD mesh used in STAR-CCM+ for computing the acoustic modes of the Helmholtz resonator.}
    \label{fig:resonator_mesh}
\end{figure}

The classical Helmholtz resonator is a very well-studied test case for modal analysis. AVSP \cite{Nicoud:2007:AIAA,Selle:2006:CNF} is a Helmholtz solver developed at CERFACS. Here we use the Helmholtz resonator geometry as described in the AVSP report \cite{Nicoud:2007:AVSP}. The resonator consists of a big cavity connected to a small open neck. In this case, the cavity is 0.1 m long, 0.08 m high, and 0.08 m deep, and the neck is a cube of size 0.02 m. There is no flow inside the resonator, the density is constant at 1.174 $kg/m^3$ and the speed of sound is also constant at 347.73 m/s. All the resonator walls are set to perfectly reflecting boundary condition, while the open end of the neck is set to zero pressure boundary condition. The CFD mesh used in STAR-CCM+ for this case is shown in \cref{fig:resonator_mesh}.

\begin{figure}
    \centering
    \begin{subfigure}[b]{0.3\textwidth}
        \centering
        \includegraphics[width=\textwidth]{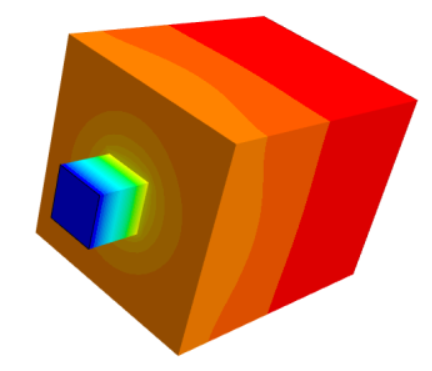}
        \caption{Mode 1}
        \label{fig:resonator_mode1}
    \end{subfigure}
    \hfill
    \begin{subfigure}[b]{0.3\textwidth}
        \centering
        \includegraphics[width=\textwidth]{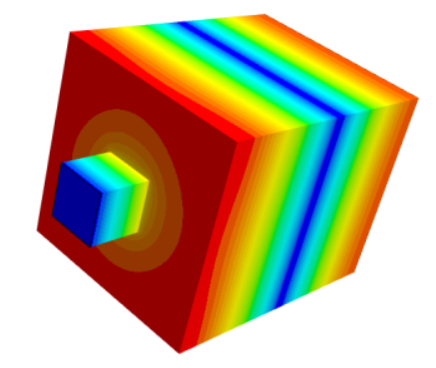}
        \caption{Mode 2}
        \label{fig:resonator_mode2}
    \end{subfigure}
    \hfill
    \begin{subfigure}[b]{0.3\textwidth}
        \centering
        \includegraphics[width=\textwidth]{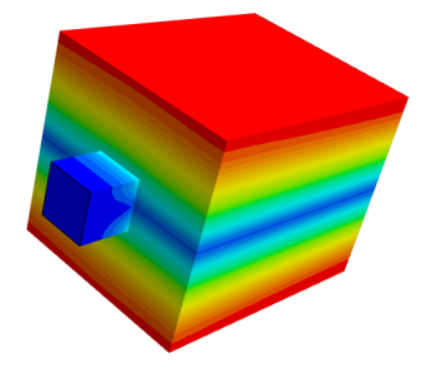}
        \caption{Mode 3}
        \label{fig:resonator_mode3}
    \end{subfigure}

    \centering
    \begin{subfigure}[b]{0.3\textwidth}
        \centering
        \includegraphics[width=\textwidth]{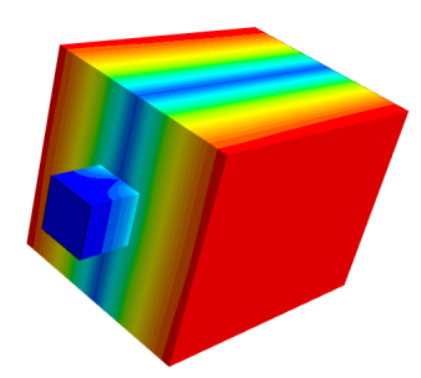}
        \caption{Mode 4}
        \label{fig:resonator_mode4}
    \end{subfigure}
    \hfill
    \begin{subfigure}[b]{0.3\textwidth}
        \centering
        \includegraphics[width=\textwidth]{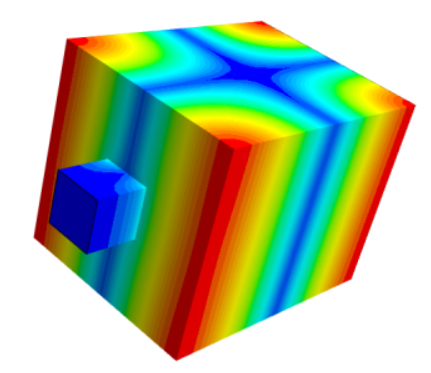}
        \caption{Mode 5}
        \label{fig:resonator_mode5}
    \end{subfigure}
    \hfill
    \begin{subfigure}[b]{0.3\textwidth}
        \centering
        \includegraphics[width=\textwidth]{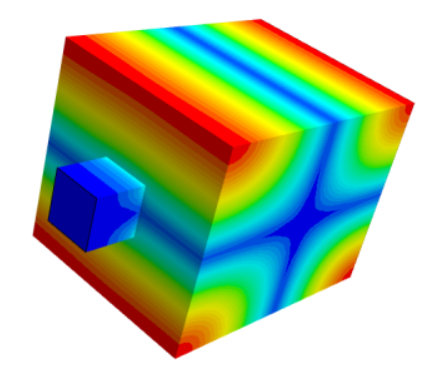}
        \caption{Mode 6}
        \label{fig:resonator_mode6}
    \end{subfigure}
     \centering
    \begin{subfigure}[b]{0.3\textwidth}
        \centering
        \includegraphics[width=\textwidth]{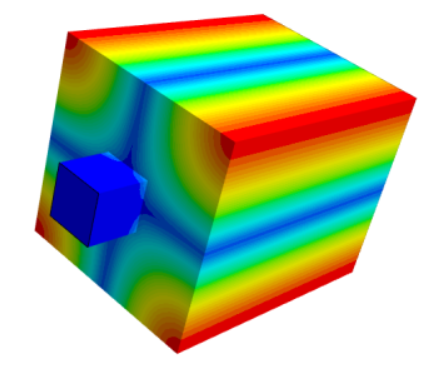}
        \caption{Mode 7}
        \label{fig:resonator_mode7}
    \end{subfigure}
    \begin{subfigure}[b]{0.3\textwidth}
        \centering
        \includegraphics[width=\textwidth]{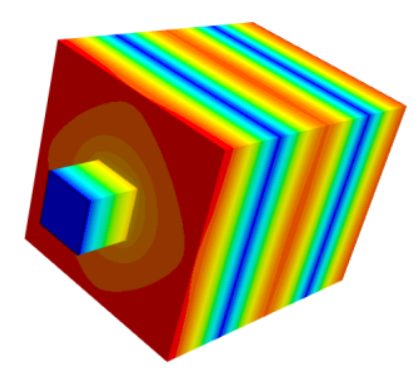}
        \caption{Mode 8}
        \label{fig:resonator_mode8}
    \end{subfigure}
    \caption{Acoustic mode shapes (normalized acoustic pressure) of the first eight modes of the Helmholtz resonator as obtained in STAR-CCM+ using SLEPc.}
    \label{fig:resonator_modes}
\end{figure}

For this particular case, the acoustic frequencies as computed by the AVSP solver \cite{Nicoud:2007:AIAA,Selle:2006:CNF} and theoretical results as reported in Table B.3 in \cite{Nicoud:2007:AVSP} have been listed in \cref{tab:resonator_results}. The acoustic frequencies computed in STAR-CCM+ using SLEPc are very similar and are also listed in \cref{tab:resonator_results}. Since this is a non-reacting case with no impedance boundary condition, the underlying eigenvalue problem is linear and is thus solved using SLEPc's \texttt{EPS} module. Also since this is a non-reacting case the growth rates are zero, i.e., we get standing waves in the resonator. The individual mode shapes are displayed in \cref{fig:resonator_modes} which match with the mode shapes reported in \cite{Nicoud:2007:AVSP}.

\begin{table}[]
    \caption{First eight acoustic mode frequencies of the Helmholtz Resonator computed in STAR-CCM+ using SLEPc, and reference results obtained from Table B.3 in \cite{Nicoud:2007:AVSP}.}
    \label{tab:resonator_results}
    \centering
    \begin{tabular}{|r|r|r|r|}
        \hline
         & \multicolumn{1}{c|}{STAR-CCM+} & \multicolumn{1}{c|}{AVSP} & \multicolumn{1}{c|}{Theoretical} \\
         & \multicolumn{1}{c|}{using SLEPc} & Table B.3 \cite{Nicoud:2007:AVSP} &  Table B.3 \cite{Nicoud:2007:AVSP} \\
         \hline
         Mode & Frequency, Hz &  Frequency, Hz & Frequency, Hz \\
         \hline
         \hline
         1 & 258 & 263 & 233 \\
         2 & 1772 & 1774 & 1781 \\
         3 & 2174 & 2176 & 2169 \\
         4 & 2174 & 2176 & 2169\\
         5 & 2785 & 2787 & 2778 \\
         6 & 2785 & 2787 & 2778 \\
         7 & 3069 & 3069 & 3068 \\
         8 & 3479 & 3479 & 3483 \\
         \hline
    \end{tabular}
\end{table}

\subsection{Rectangular Cavity}\label{sec:cavity}
To validate the Helmholtz equation solution with impedance boundary condition, we consider the rectangular cavity case described in \cite{Nicoud:2007:AIAA}. The 2D cavity has a length, $L = 0.5$ m and height $h = 0.1$ m with a constant speed of sound $c_0 = 450$ m/s. The left, top and bottom boundaries of the cavity are set to perfectly reflecting boundary condition, while the right boundary is set to the impedance boundary condition with a specified constant impedance value of $Z = a + \iota b$.

For this particular case, the Helmholtz equation exhibits 1D longitudinal modes which can be solved analytically and the acoustic frequencies are given as
\begin{equation}
\label{eq:cavity_fm}
    f_m = m \frac{c_o}{2L} + \frac{c_0}{2 \pi L} \arctan\left(\frac{-\iota}{Z}\right),
\end{equation}
where $m$ is the mode number.

For the case where the impedance is purely reactive (imaginary) impedance, $Z = \iota b$, we get purely real frequencies
\begin{equation}
\label{eq:fm_b}
    f_m = m \frac{c_o}{2L} + \frac{c_0}{2 \pi L} \arctan\left(\frac{-1}{b}\right),
\end{equation}
and for the case where we have purely resistive (real) impedance, $Z = a$, we get the following complex frequencies
\begin{equation}
\label{eq:fm_a}
    f_m = m \frac{c_o}{2L} - \iota \frac{c_0}{4 \pi L} \ln\left(\frac{a+1}{a-1}\right).
\end{equation}

\begin{figure}
    \centering
    \includegraphics[width=400pt]{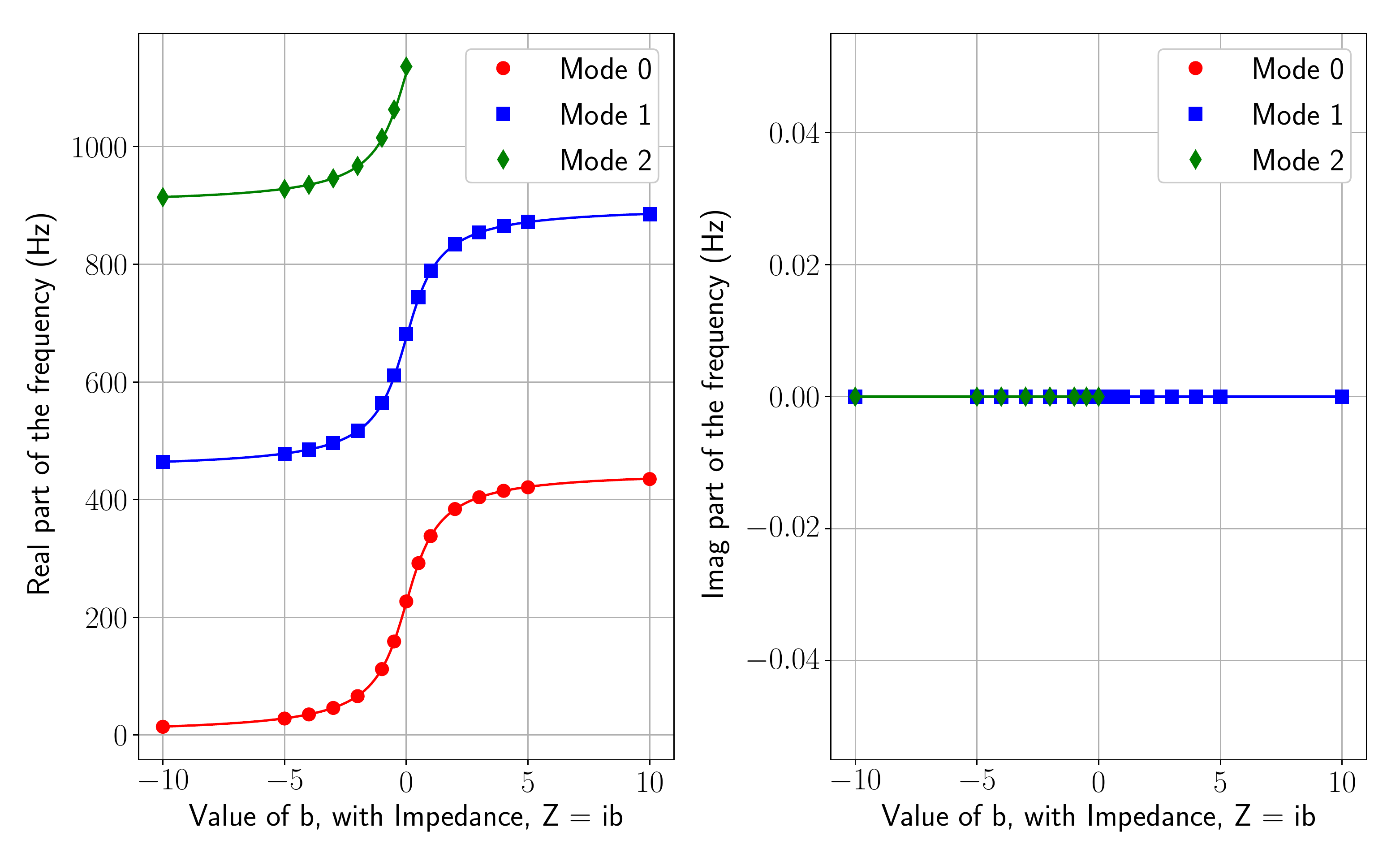}
    \caption{Acoustic mode frequencies as a function of purely reactive (imaginary) impedance boundary condition, $Z = \iota b$. The dashed line is the analytical solution given by \cref{eq:fm_b} and the symbols are the modes computed in STAR-CCM+ using SLEPc.}
    \label{fig:Z_b_impedance}
\end{figure}

\begin{figure}
    \centering
    \includegraphics[width=400pt]{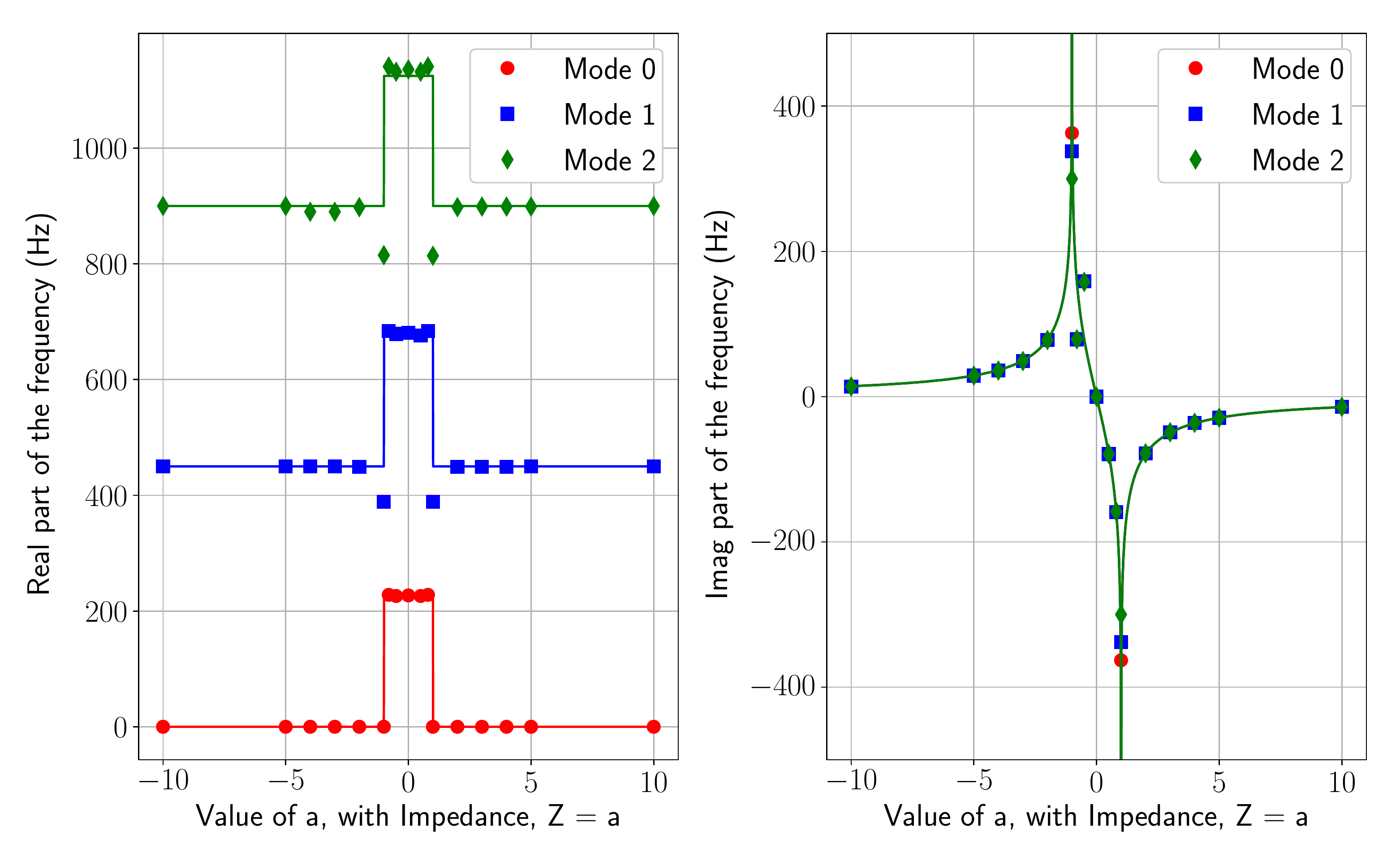}
    \caption{Acoustic mode frequencies as a function of purely resistive (real) impedance boundary condition, $Z = a$. The dashed line is the analytical solution given by \cref{eq:fm_a} and the symbols are the modes computed in STAR-CCM+ using SLEPc.}
    \label{fig:Z_a_impedance}
\end{figure}

In \cite{Nicoud:2007:AIAA}, the acoustic frequencies for this case are computed numerically using the AVSP solver for purely reactive and purely resistive impedance values (results are presented in Figure 2 and Figure 3 in \cite{Nicoud:2007:AIAA}). We repeated the same set of tests using our solver with SLEPc. In this case since the underlying problem is quadratic, we use SLEPc's \texttt{PEP} package to compute the eigenvalues. The results computed using SLEPc are presented in \cref{fig:Z_b_impedance,fig:Z_a_impedance}. We see that the numerical frequencies computed using SLEPc match the analytical results given by \cref{eq:fm_b,eq:fm_a}.

For the purely reactive impedance, $Z = \iota b$, as the impedance value approaches zero, $|b|\rightarrow 0$, the $1/b$ term in \cref{eq:fm_b} makes the numerical solution ill-conditioned. We sometimes get a few additional transverse modes due to numerical noise, but they can be easily filtered out by looking at the acoustic mode shapes. Similarly for the case of purely resistive impedance, $Z = a$, we get some numerical instability around $a = 1$ where the acoustic frequency is not well defined. For this reason, we see some oscillations at $Z = \pm 1$ in  \cref{fig:Z_a_impedance}. Also for the case when $b = 0$ or $a = 0$ we get $Z = 0$, which reduces to the zero acoustic pressure boundary condition. So when $Z = 0$ we simply fall back to the zero acoustic pressure boundary condition and the eigenvalue problem becomes linear.

\subsection{1D Planar Flame}\label{sec:planar_flame}
\begin{figure}
    \centering
    \includegraphics[width=400pt]{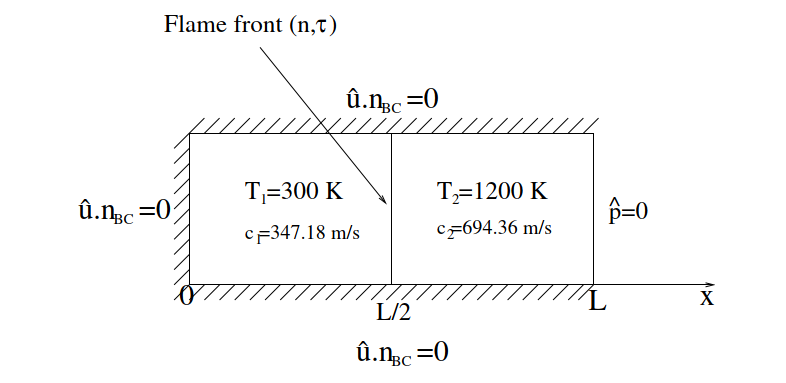}
    \caption{Schematic of the 1D planar duct flame extracted from \cite{Nicoud:2007:AIAA}.}
    \label{fig:1dflame_sketch}
\end{figure}

Here we consider a simple planar flame case in a 1D duct as shown in \cref{fig:1dflame_sketch} which is extracted from \cite{Nicoud:2007:AIAA}. A thin flame in the middle of the duct separates two regions of cold fresh unburnt gases to the left and the hot burnt gases to the right. For this particular test case, assuming an infinitely thin flame and jump conditions across the flame allows the acoustic modes to be computed analytically as is shown in \cite{Nicoud:2007:AIAA,Kaufmann:2002:CNF}.

\begin{figure}
    \centering
    \includegraphics[width=400pt]{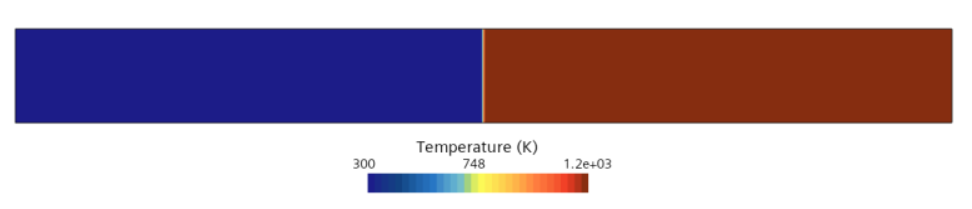}
    \caption{Temperature field across a thin flame stabilized for simulating the 1D planar flame in STAR-CCM+.}
    \label{fig:1dflame_temp}
\end{figure}

\begin{figure}
    \centering
    \begin{subfigure}[b]{0.45\textwidth}
        \centering
        \includegraphics[width=\textwidth]{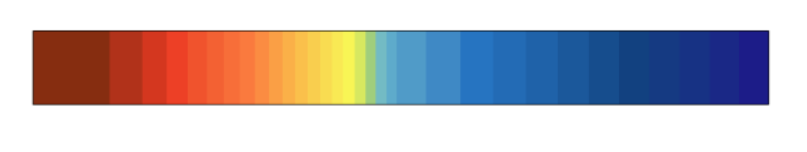}
        \caption{Mode 1}
    \end{subfigure}
    \hfill
    \begin{subfigure}[b]{0.45\textwidth}
        \centering
        \includegraphics[width=\textwidth]{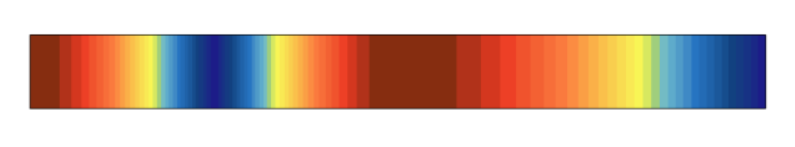}
        \caption{Mode 2}
    \end{subfigure}

    \centering
    \begin{subfigure}[b]{0.45\textwidth}
        \centering
        \includegraphics[width=\textwidth]{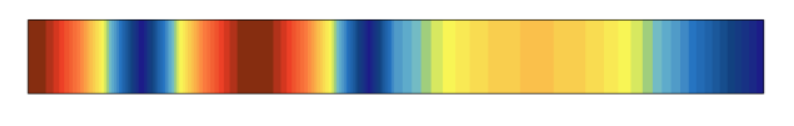}
        \caption{Mode 3}
    \end{subfigure}
    \hfill
    \begin{subfigure}[b]{0.45\textwidth}
        \centering
        \includegraphics[width=\textwidth]{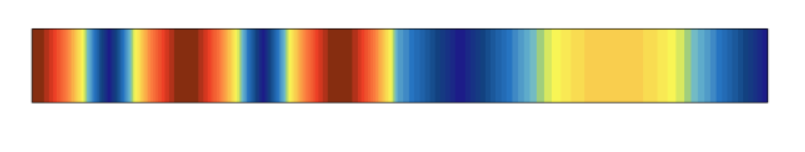}
        \caption{Mode 4}
    \end{subfigure}
    \caption{First four acoustic mode shapes of the 1D planar duct flame case obtained in STAR-CCM+.}
    \label{fig:1dflame_modes}
\end{figure}

The acoustic modes of this planar flame configuration are computed numerically in STAR-CCM+. In the simulation, a thin 1D flame heat source is imposed along with a temperature distribution as shown in \cref{fig:1dflame_temp} to match the conditions described in \cref{fig:1dflame_sketch}. A relatively coarse mesh containing 1000 cells is used. For the Helmholtz solver, the inlet boundary is set to perfectly reflecting condition, while the outlet boundary is set to zero acoustic pressure. The simplified $n-\tau$ model as described in \cref{sec:n-tau} is used to compute the heat source term with a constant $\tau = 10^{-4}$ s and an interaction index coefficient $\eta = 5$ (as used in \cite{Nicoud:2007:AIAA}). Due to the presence of the heat source term, the underlying eigenvalue problem is nonlinear and so we use SLEPc's \texttt{NEP} module to solve this problem. The first four acoustic modes computed using SLEPc are listed in \cref{tab:1dflame_results}, and the mode shapes are displayed in \cref{fig:1dflame_modes}. These results are consistent with the results reported in \cite{Nicoud:2007:AIAA}. The first and fourth modes are stable, the second mode is marginally stable, and the third mode is unstable. The results are off from the analytical results obtained for the infinitely thin flame, as the numerical flame on the coarse grid used in this simulation is thicker, however, the trend of the growth rate vs. mode number is consistent. Localized mesh refinement near the flame as is done in \cite{Nicoud:2007:AIAA} should yield better results.

\begin{table}[]
    \caption{First four acoustic mode frequencies and growth rates for the 1D planar duct flame case as obtained by the Helmholtz solver implemented in STAR-CCM+ using SLEPc. The first and fourth modes that have a negative growth rate are stable, the second one is marginally stable, and the third mode is unstable with a positive growth rate. Reference analytical results as obtained in \cite{Nicoud:2007:AIAA} for an infinitely thin flame are also provided.}
    \label{tab:1dflame_results}
    \centering
    \begin{tabular}{|r||r|r||r|r|}
        \hline
        & \multicolumn{2}{c||}{STAR-CCM+} & \multicolumn{2}{c|}{Analytical Eq.(56) \cite{Nicoud:2007:AIAA}} \\
        \hline
        \hline
         Acoustic & Frequency & Growth Rate & Frequency & Growth Rate \\
         Mode & $f_r$ (Hz) & $\omega_i$ (rad/s) & $f_r$ (Hz) & $\omega_i$ (rad/s) \\
         \hline
         \hline
         1 & 265.2  & -7.08 & 159.6 & -32.91 \\
         2 & 694.7 &  0.10 & 694.4 & 0.0 \\
         3 & 1124.8 & 29.94 & 1227.3 & 261.67 \\
         4 & 1659.8 & -42.65 & 1546.6 & -336.76 \\
         \hline
    \end{tabular}
\end{table}

\subsection{BH3P Burner}\label{sec:bh3p}
\begin{figure}
    \centering
    \includegraphics[width=400pt]{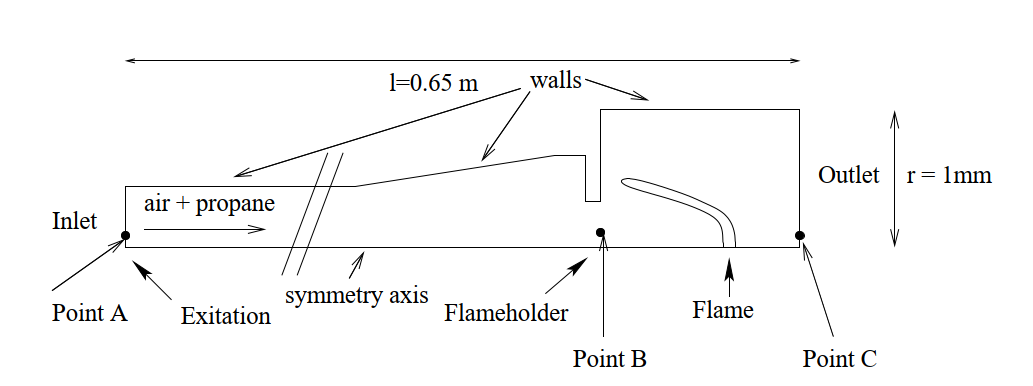}
    \caption{Schematic (not to scale) of the BH3P burner as shown in \cite{Kaufmann:2002:CNF}.}
    \label{fig:bh3p_sketch}
\end{figure}

The BH3P burner is described in the works of \cite{Kaufmann:2002:CNF, Truffin:2005:CNF}. It consists of a ducted premixed propane-air flame which is stabilized downstream using a perforated plate as shown in the burner geometry \cref{fig:bh3p_sketch}.

\begin{figure}
    \centering
    \includegraphics[width=400pt]{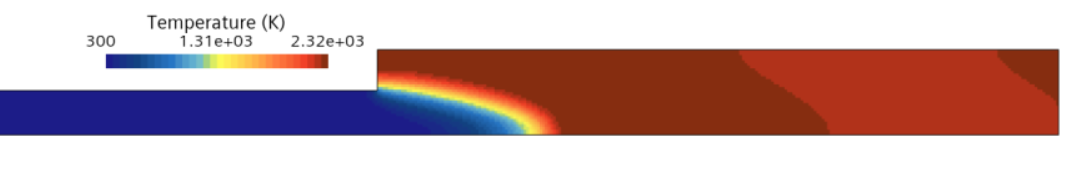}
    \includegraphics[width=400pt]{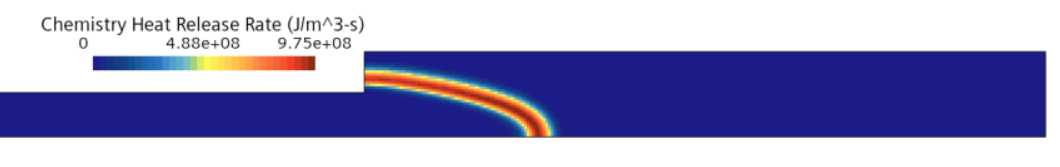}
    \caption{Temperature field across the converged flame (top) and the chemistry heat release rate (bottom) obtained in STAR-CCM+ simulating the BH3P burner case. The scene is zoomed in near the burner, the upstream channel is not shown.}
    \label{fig:bh3p_temp}
\end{figure}

We first simulate this flame configuration in STAR-CCM+ using a single-step reaction using the Eddy Breakup Model to get a stable flame in approximately the same location as described in \cite{Kaufmann:2002:CNF, Truffin:2005:CNF}. The temperature field around the flame is shown in \cref{fig:bh3p_temp}. We use our acoustic Helmholtz solver to compute the acoustic frequencies and growth rates using the converged solution obtained from the CFD simulation of the BH3P burner. We use our simplified $n-\tau$ model to compute the heat release rate with a constant $\tau = 5 \times 10^{-4}$ s and $\eta = 4$ (similar to values reported in Fig.~7 in \cite{Truffin:2005:CNF}). The inlet and side walls are set to perfectly reflecting, while the outlet is set to zero acoustic pressure boundary condition.

\begin{figure}
    \centering
    \begin{subfigure}[b]{0.45\textwidth}
        \centering
        \includegraphics[width=\textwidth]{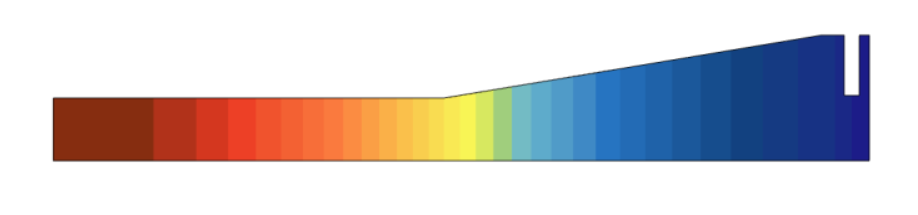}
        \caption{Mode 1}
        \label{fig:bh3p_mode1}
    \end{subfigure}
    \hfill
    \begin{subfigure}[b]{0.45\textwidth}
        \centering
        \includegraphics[width=\textwidth]{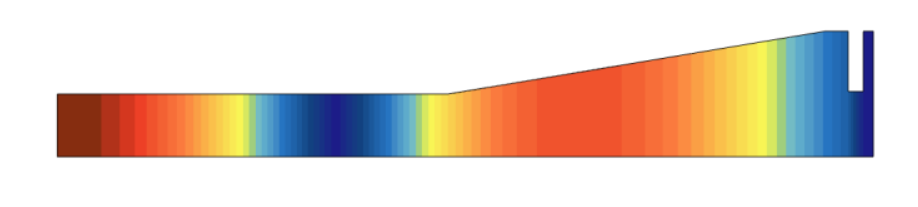}
        \caption{Mode 2}
        \label{fig:bh3p_mode2}
    \end{subfigure}

    \centering
    \begin{subfigure}[b]{0.45\textwidth}
        \centering
        \includegraphics[width=\textwidth]{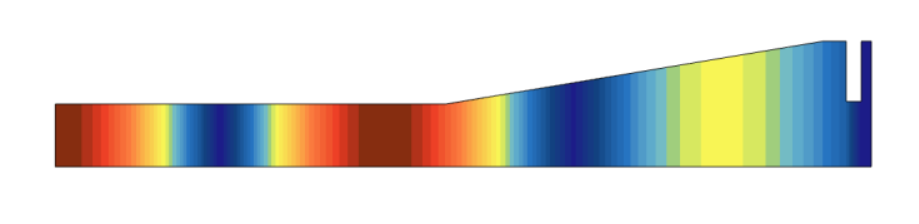}
        \caption{Mode 3}
        \label{fig:bh3p_mode3}
    \end{subfigure}
    \hfill
    \begin{subfigure}[b]{0.45\textwidth}
        \centering
        \includegraphics[width=\textwidth]{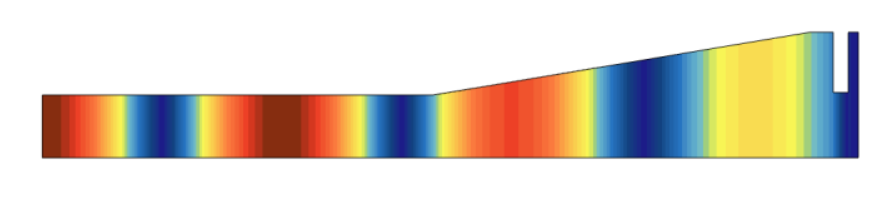}
        \caption{Mode 4}
        \label{fig:bh3p_mode4}
    \end{subfigure}

    \centering
    \begin{subfigure}[b]{0.45\textwidth}
        \centering
        \includegraphics[width=\textwidth]{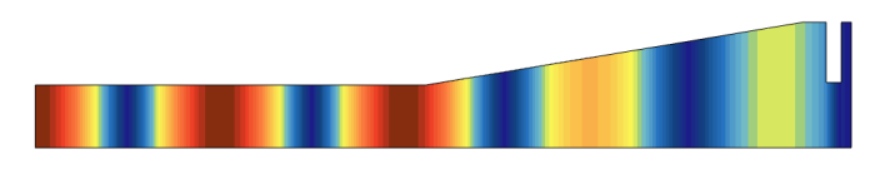}
        \caption{Mode 5}
        \label{fig:bh3p_mode5}
    \end{subfigure}
    \hfill
    \begin{subfigure}[b]{0.45\textwidth}
        \centering
        \includegraphics[width=\textwidth]{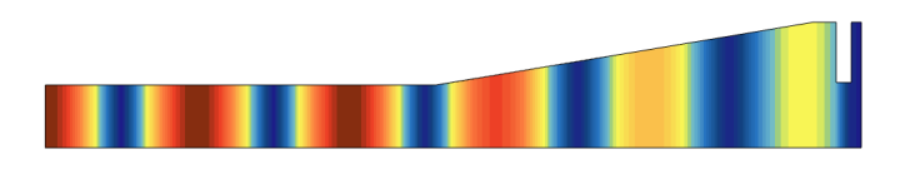}
        \caption{Mode 6}
        \label{fig:bh3p_mode6}
    \end{subfigure}
    \caption{First six acoustic mode shapes of the BH3P burner obtained in STAR-CCM+.}
    \label{fig:bh3p_modes}
\end{figure}

The first six computed modes are listed in \cref{tab:bh3p_results} and the corresponding mode shapes are shown in \cref{fig:bh3p_modes}. Also presented in the same table are results obtained from \cite{Kaufmann:2002:CNF}, and as we can see, we get very similar results. The fifth and sixth modes are unstable, while the first four modes are stable.

\begin{table}[]
    \caption{First six acoustic mode frequencies and growth rates for the BH3P burner flame as obtained by the Helmholtz solver implemented in STAR-CCM+ using SLEPc and reference results obtained from \cite{Kaufmann:2002:CNF}.}
    \label{tab:bh3p_results}
    \centering
    \begin{tabular}{|r||r|r||r|r|}
        \hline
        & \multicolumn{2}{c||}{STAR-CCM+} & \multicolumn{2}{c|}{Results from \cite{Kaufmann:2002:CNF}} \\
        \hline
        \hline
         Acoustic & Frequency & Growth Rate & Frequency & Growth Rate \\
         Mode & $f_r$ (Hz) & $\omega_i$ (rad/s) & $f_r$ (Hz) & $\omega_i$ (rad/s) \\
         \hline
         \hline
         1 & 168.7  & -1.98 & 159.5 & -1.26 \\
         2 & 382.2 & -25.58 & 417.5 & -28.90 \\
         3 & 644.6 & -26.37 & 607.5 & -35.18 \\
         4 & 887.5 &  -14.32 & 881.7 & -13.19 \\
         5 & 1154.2 & 17.51 & 1118.6 & 21.36 \\
         6 & 1399.6 & 36.25 & 1354.0 & 38.33 \\
         \hline
    \end{tabular}
\end{table}

\subsection{EM2C Burner}\label{sec:em2c}
\begin{figure}
    \centering
    \includegraphics[width=400pt]{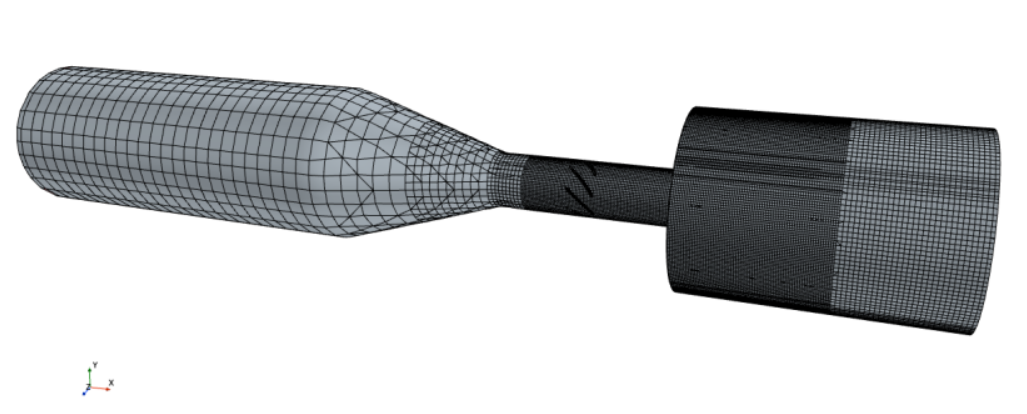}
    \caption{The computational mesh used for the EM2C burner case.}
    \label{fig:em2c_mesh}
\end{figure}

\begin{figure}
    \centering
    \includegraphics[width=400pt]{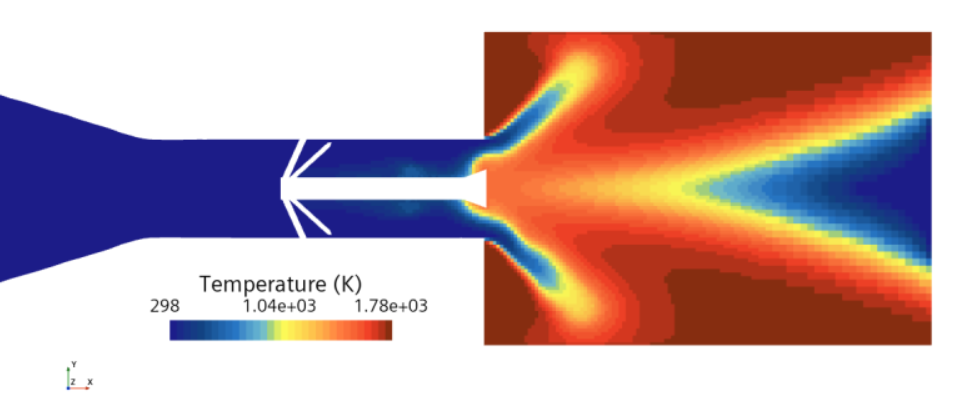}
    \includegraphics[width=400pt]{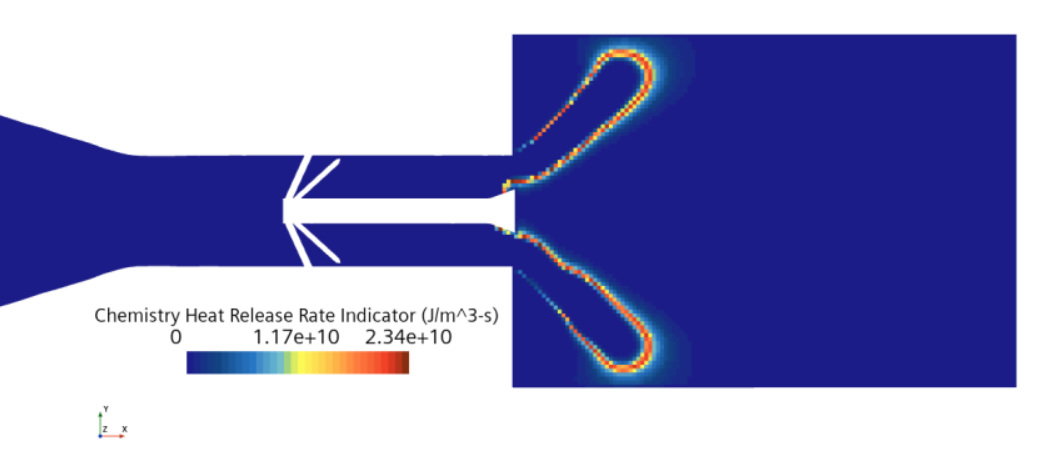}
    \caption{The temperature field (top) and chemistry heat release rate (bottom) obtained through reacting flow simulation of the EM2C burner. The scene is zoomed in near the flame position, the rest of the channel is not shown.}
    \label{fig:em2c_temp}
\end{figure}

\begin{figure}
    \centering
    \begin{subfigure}[b]{0.45\textwidth}
        \centering
        \includegraphics[width=\textwidth]{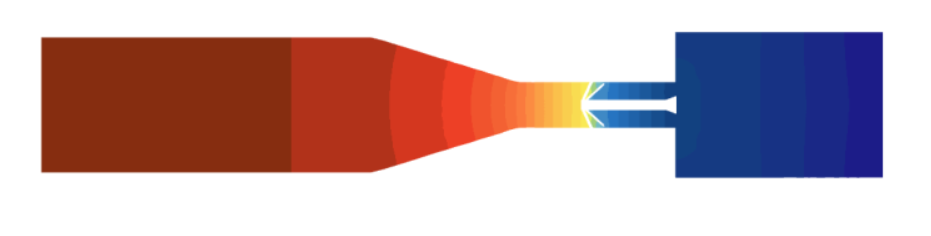}
        \caption{Mode 1}
        \label{fig:em2c_mode1}
    \end{subfigure}
    \hfill
    \begin{subfigure}[b]{0.45\textwidth}
        \centering
        \includegraphics[width=\textwidth]{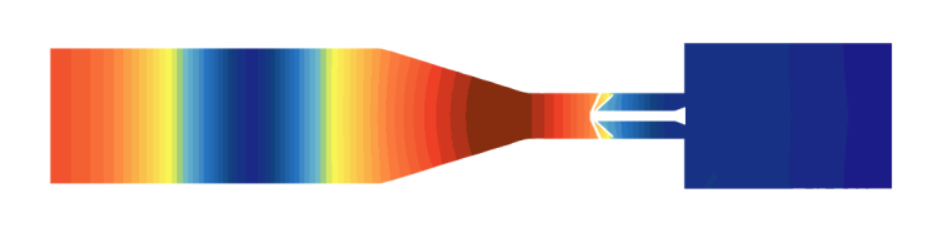}
        \caption{Mode 2}
        \label{fig:em2c_mode2}
    \end{subfigure}

    \centering
    \begin{subfigure}[b]{0.45\textwidth}
        \centering
        \includegraphics[width=\textwidth]{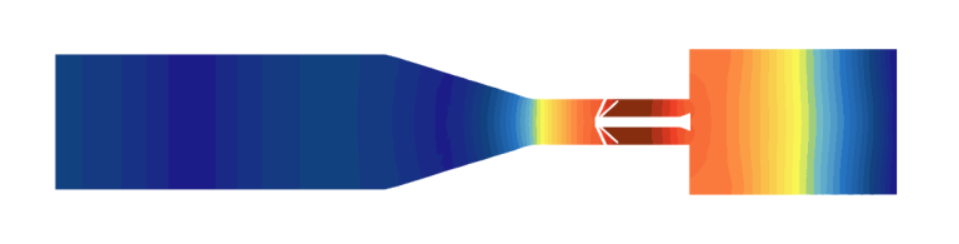}
        \caption{Mode 3}
        \label{fig:em2c_mode3}
    \end{subfigure}
    \hfill
    \begin{subfigure}[b]{0.45\textwidth}
        \centering
        \includegraphics[width=\textwidth]{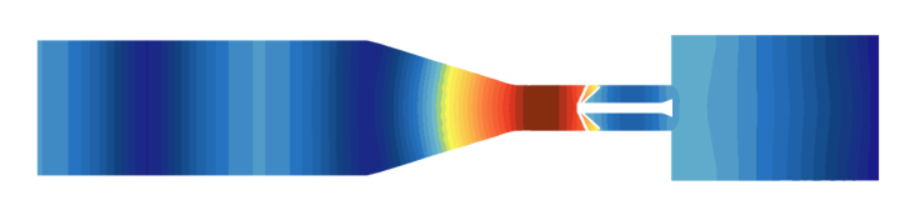}
        \caption{Mode 4}
        \label{fig:em2c_mode4}
    \end{subfigure}
    \caption{First four acoustic mode shapes of the EM2C burner.}
    \label{fig:em2c_modes}
\end{figure}

Here we consider the turbulent swirled combustor designed and studied at the EM2C Laboratory which has been extensively studied in many previous works \cite{Palies:2010:CNF, Silva:2013:CNF, Falco:2021:Thesis}. The combustor geometry comprises a variable length upstream manifold, an axisymmetric duct and a variable length cylindrical combustion chamber. The upstream manifold and combustion cylinder lengths are varied to study different flame configurations as described in \cite{Silva:2013:CNF}. We pick the CO5 configuration to study in this case. We simulate this EM2C combustor (CO5 conditions) in STAR-CCM+ using the Flamelet Generated Manifold (FGM) reacting model. The mesh used for simulation in shown in \cref{fig:em2c_mesh}. The chemistry heat release rate and the temperature field around the converged flame is shown in \cref{fig:em2c_temp}. The Helmholtz solver along with the $n-\tau$ model is used to compute the acoustic modes using the converged solution. In \cite{Palies:2010:CNF}, the Gain $G$ (similar to $\eta$) and phase $\phi$ (where $\phi = \omega \tau$) for this flame are measured experimentally as a function of acoustic perturbation frequency and this data is shown in Fig.3 in \cite{Silva:2013:CNF}. Since in our simplified $n-\tau$ model, we use a constant $\eta$ and $\tau$, in this study we fix the interaction index coefficient $\eta = 1$ and vary the time lag $\tau$ to get the first acoustic frequency close to the reported value of $f_{1st} = 118.4$ Hz in \cite{Silva:2013:CNF}. Based on some testing we have used $\tau = 7.3 \times 10^{-4}$ in this study. The inlet and walls are set to perfectly reflecting, and the outlet is set to zero acoustic pressure boundary condition. The first four computed acoustic frequencies and growth rates are listed in \cref{tab:em2c_results} and the corresponding mode shapes are shown in \cref{fig:em2c_modes}. Based on the results, we find that the first and fourth modes are stable, while the second and third modes are unstable. Consistent to these results, in \cite{Silva:2013:CNF} the first mode of CO5 configuration is found to be stable.

\begin{table}[]
    \caption{First four acoustic mode frequencies and growth rates of the EM2C burner case.}
    \label{tab:em2c_results}
    \centering
    \begin{tabular}{|r|r|r|}
        \hline
         Acoustic & Frequency & Growth Rate \\
         Mode & $f_r$ (Hz) & $\omega_i$ (rad/s) \\
         \hline
         \hline
         1 & 116.6 & -37.73 \\
         2 & 907.5 & 15.31 \\
         3 & 1334.9 & 336.01 \\
         4 & 1648.0 &  -173.22 \\
         \hline
    \end{tabular}
\end{table}

\section{Comparison of different algorithms}
As described in the earlier sections there are different ways of solving the nonlinear eigenvalue problem. In many previous implementations of the Helmholtz solver, the linearized iterative algorithm as described in \cref{sec:lin_iter} is used. In SLEPc's \texttt{NEP} module many different algorithms are implemented and in this work we use the NLEIGS algorithm \cite{Guttel:2014:NCF}. In this section we compare these two algorithms for the BH3P and EM2C burner cases described in the previous sections.

\subsection{BH3P Burner}
Here we use the BH3P burner case as described in \cref{sec:bh3p}. For both the iterative and NLEIGS algorithms we first need to compute the acoustic modes without the reaction source. The cold flow acoustic modes for the BH3P case are listed in \cref{tab:bh3p_comparison}. In the NLEIGS algorithm, we need to provide a region in the complex frequency space to compute the acoustic modes. We estimate this region size based on the minimum and maximum frequency computed for the cold flow. Given this region, NLEIGS algorithm can compute the required number of eigenvalues within that region. So using this procedure we compute the acoustic modes of the full nonlinear eigenvalue problem using SLEPc. The computed modes are listed in \cref{tab:bh3p_comparison}. On the other hand, the iterative algorithm as described in \cref{sec:lin_iter} solves the linearized eigenvalue problem iteratively for each mode starting from the cold flow frequency. The results computed using this iterative algorithm are also listed in the same table.

As we see from the results listed in \cref{tab:bh3p_comparison}, for this BH3P case both SLEPc's NLEIGS algorithm and the iterative algorithm yield very similar acoustic modes. However, SLEPc's NLEIGS is almost $4\times$ faster compared to the iterative algorithm as we need to repeatedly solve the linearized system at each iteration and typically each mode takes 2-3 iterations to converge.

\begin{table}[]
    \caption{Comparison of results of the BH3P burner case computed using SLEPc's NLEIGS algorithm against an iterative algorithm.}
    \label{tab:bh3p_comparison}
    \centering
    \begin{tabular}{|r||r|r||r|r||r|r|}
        \hline
         {} & \multicolumn{2}{c||}{Cold Flow} & \multicolumn{2}{c||}{SLEPc's NLEIGS} & \multicolumn{2}{c|}{Iterative} \\
         \hline
         \hline
         Acoustic & Frequency & Growth Rate & Frequency & Growth Rate & Frequency & Growth Rate \\
         Mode &  $f_r$ (Hz) & $\omega_i$ (rad/s) & $f_r$ (Hz) & $\omega_i$ (rad/s) & $f_r$ (Hz) & $\omega_i$ (rad/s) \\
         \hline
         \hline
         1 & 169.4 & 0.0 & 168.7 & -2.79 & 168.7 & -2.77 \\
         2 & 383.8 & 0.0 & 382.2 & -26.10 & 382.2 & -26.08\\
         3 & 642.5 & 0.0 & 644.5 & -26.23 & 644.5 & -26.23\\
         4 & 881.5 & 0.0 & 887.5 & -14.19 & 887.5 & -14.19 \\
         5 & 1149.2 & 0.0 & 1154.2 & 17.54 & 1154.3 & 17.54\\
         6 & 1398.0 & 0.0 & 1399.6 & 36.23 & 1399.6 & 36.23 \\
         \hline
         \hline
         Runtime & \multicolumn{2}{c||}{5 s} & \multicolumn{2}{c||}{15 s} & \multicolumn{2}{c|}{58 s} \\
         \hline
    \end{tabular}
\end{table}

\subsection{EM2C Burner}
As in the previous section, here we present results for the EM2C burner case described in \cref{sec:em2c}. The acoustic modes for the cold flow, and for the full reacting case computed using SLEPc's NLEIGS algorithm and the iterative algorithm are listed in \cref{tab:em2c_comparison}. In this case we notice the NLEIGS is able to find the four modes, however the iterative algorithm struggles to converge for the 3rd mode. One likely reason is that starting from the 3rd mode's cold flow frequency of $f_r = 1564.8$ Hz the reactive modes lie very close to the left (1334.2 Hz) and right (1640.2 Hz) of that frequency. So when using the iterative algorithm, the solution keeps oscillating between these modes and is unable to converge. Also, the 4th mode obtained using the iterative algorithm is close to the 5th mode computed using NLEIGS (not listed in the table). So in summary, we are unable to compute the 3rd and 4th modes using the iterative algorithm. Moreover, the iterative algorithm takes almost $6\times$ the time compared to NLEIGS to compute the modes.

In \cref{tab:em2c_iter_conv}, we show the convergence of individual modes using the iterative algorithm. Typically most modes converge in 2-3 modes as we can see with the 1st and 2nd modes. However, the 3rd mode keeps oscillating and does not converge even after 5 iterations. The 4th mode converges but to the 5th mode computed using SLEPc. The iterative algorithm was rerun for 10 iterations, but still the 3rd mode did not converge.

\subsection{Contour integration and other algorithms}
The iterative algorithm can sometimes struggle to converge and miss modes as has also been reported in \cite{Buschmann:2020:JEGTP}. In \cite{Buschmann:2020:JEGTP}, an alternative contour-integration based method is proposed based on the Beyn's algorithm \cite{Beyn:2012:LAA} to compute the eigenvalues. In \cite{Mohamed:2020:EABM}, a rational approximation method using Raleigh-Ritz algorithm is used for solving the acoustic nonlinear problem. Some of the these contour-integration based algorithms are also implemented in SLEPc under the Contour Integral Spectrum Slicing (CISS) method \cite{Maeda:2016:CIS} along with Hankel and Raleigh-Ritz extraction methods \cite{Yokota:2013:PMN,Imakura:2018:BSC}. In our preliminary testing we found NLEIGS method to perform as good or better than the CISS based methods, so in this study we have exclusively used the NLEIGS method.

From our previous experience in the context of other problems, contour integral methods are not generally competitive with respect to Krylov methods such as TOAR or NLEIGS. On one hand, contour integration requires a factorization at each integration point (which can be 32 or even more depending on the problem), while NLEIGS computes just one factorization (and still this is the most computationally expensive step of the overall method). Furthermore, in contour integral methods the factorized matrix is applied to a block of vectors, that is, multiple right-hand-side solves are required with a number of vectors that must be larger than the number of eigenvalues inside the contour. This is also an increased cost with respect to NLEIGS, which operates as a single-vector recurrence. The accuracy of contour integral methods depends on many parameters that may be difficult to choose, for instance the number of integration points, the size of the subspace, or a threshold tolerance when estimating the number of eigenvalues inside the contour. If not chosen appropriately, the solver may require one or more refinement steps whenever the accuracy is poor with respect to the user tolerance. This multiplies the overall cost. As a strong point, with contour integral methods there is the guarantee that all eigenvalues inside the region are returned, while with NLEIGS the user specifies how many eigenvalues are needed, without knowing if the actual number is larger or smaller. Also, contour integration admits more parallelism and therefore it is possible to implement it in a scalable way, so that the higher cost is compensated if more processors are used.

\begin{table}[]
    \caption{Comparison of acoustic modes of the EM2C burner case computed using SLEPc's NLEIGS algorithm against an iterative algorithm.}
    \label{tab:em2c_comparison}
    \centering
    \begin{tabular}{|r||r|r||r|r||r|r|}
        \hline
         {} & \multicolumn{2}{c||}{Cold Flow} & \multicolumn{2}{c||}{SLEPc's NLEIGS} & \multicolumn{2}{c|}{Iterative} \\
         \hline
         \hline
         Acoustic & Frequency & Growth Rate & Frequency & Growth Rate & Frequency & Growth Rate \\
         Mode &  $f_r$ (Hz) & $\omega_i$ (rad/s) & $f_r$ (Hz) & $\omega_i$ (rad/s) & $f_r$ (Hz) & $\omega_i$ (rad/s) \\
         \hline
         \hline
         1 & 120.1 & 0.0 & 115.4 & -36.93 & 115.4 & -36.93\\
         2 & 903.1 & 0.0 & 906.3 & 14.81 & 906.3 & 14.81 \\
         3 & 1564.8 & 0.0 & 1334.2 & 341.98 & 948.8 & -10407.97\\
         4 & 1856.6 & 0.0 & 1640.2 & -173.60 & 2150.9 & -166.61 \\
         \hline
         \hline
         Runtime & \multicolumn{2}{c||}{270 s} & \multicolumn{2}{c||}{1210 s} & \multicolumn{2}{c|}{6820 s} \\
         \hline
    \end{tabular}
\end{table}

\begin{table}[]
    \caption{Convergence of acoustic modes of EM2C burner computed using the iterative algorithm.}
    \label{tab:em2c_iter_conv}
    \centering
    \begin{tabular}{|r||r|r||r|r||r|r||r|r|}
        \hline
        & \multicolumn{2}{c||}{Mode 1} & \multicolumn{2}{c||}{Mode 2} & \multicolumn{2}{c||}{Mode 3} & \multicolumn{2}{c|}{Mode 4} \\
        \hline
        \hline
        Iteration & $f_r$ & $\omega_i$ & $f_r$ & $\omega_i$ & $f_r$ & $\omega_i$ & $f_r$ & $\omega_i$\\
        \hline
        \hline
        0 & 120.1 & 0.0 & 903.1 & 0.0 & 1564.8 & 0.0 & 1856.6 & 0.0 \\
        1 & 115.8 & -35.54 & 906.3 & 14.61 & 1429.9 & -2185.64 & 2003.0 & -2860.90 \\
        2 & 115.5 & -36.89 & 906.3 & 14.81 & 915.9 & -9.94 & 2131.8 & 4.71\\
        3 & 115.4 & -36.93 & 906.3 & 14.81 & 1620.1 & 114.63 & 2157.2 & -214.37\\
        4 & 115.4 & -36.93 & 906.3 & 14.81 & 1551.7 & -2335.68 & 2148.2 & -154.24\\
        5 & 115.4 & -36.93 & 906.3 & 14.81 &  948.8 & -10407.97 & 2150.9 & -166.61 \\
        \hline
    \end{tabular}
\end{table}

\section{Extension to non-constant $n$ and $\tau$}
In this work, we have used the $n-\tau$ model with a constant interaction index $\eta$ and time lag $\tau$ to model the heat release rate as described in \cref{sec:const-tau}. Here we describe some ideas to extend this to variable interaction index and time lag.

One simple extension of the $n-\tau$ model to accommodate multiple sources of disturbances is to use the Distributed Time Delay (DTD) method as proposed in \cite{Poilifke:2020:PECS} which expresses the heat release source as a linear combination of many individual $n-\tau$ models as follows
\begin{equation}
    D(\omega) = \sum_i e^{\iota \omega \tau_i} S_i.
\end{equation}

The DTD method allows us to use multiple distributed $n-\tau$ models each with its own interaction index $\eta_i$ and time lag $\tau_i$. This linear combination DTD model can be easily incorporated into SLEPc's NLEIGS algorithm following \cref{eq:split} to compute the acoustic modes of the combined DTD system.

\section{Scaling Studies}
Here we present results of some scaling studies we performed to study the performance (runtime) and memory consumption of our acoustic modal solver as a function of mesh size (number of cells) and in parallel using different number of processes. We consider here two test cases: the Helmholtz Resonator and the BH3P burner. The mesh scaling studies are performed in serial (single process) on an Intel Xeon Gold 6130 CPU 2.10 GHz processor with a maximum memory of 120 GB. And the parallel scaling studies are performed on a cluster where each node has 40 processors of Intel Xeon Gold 6248 CPU 2.50 GHz with a maximum memory of 180 GB per node.

\subsection{Helmholtz Resonator}
\begin{figure}
    \centering
    \includegraphics[width=400pt]{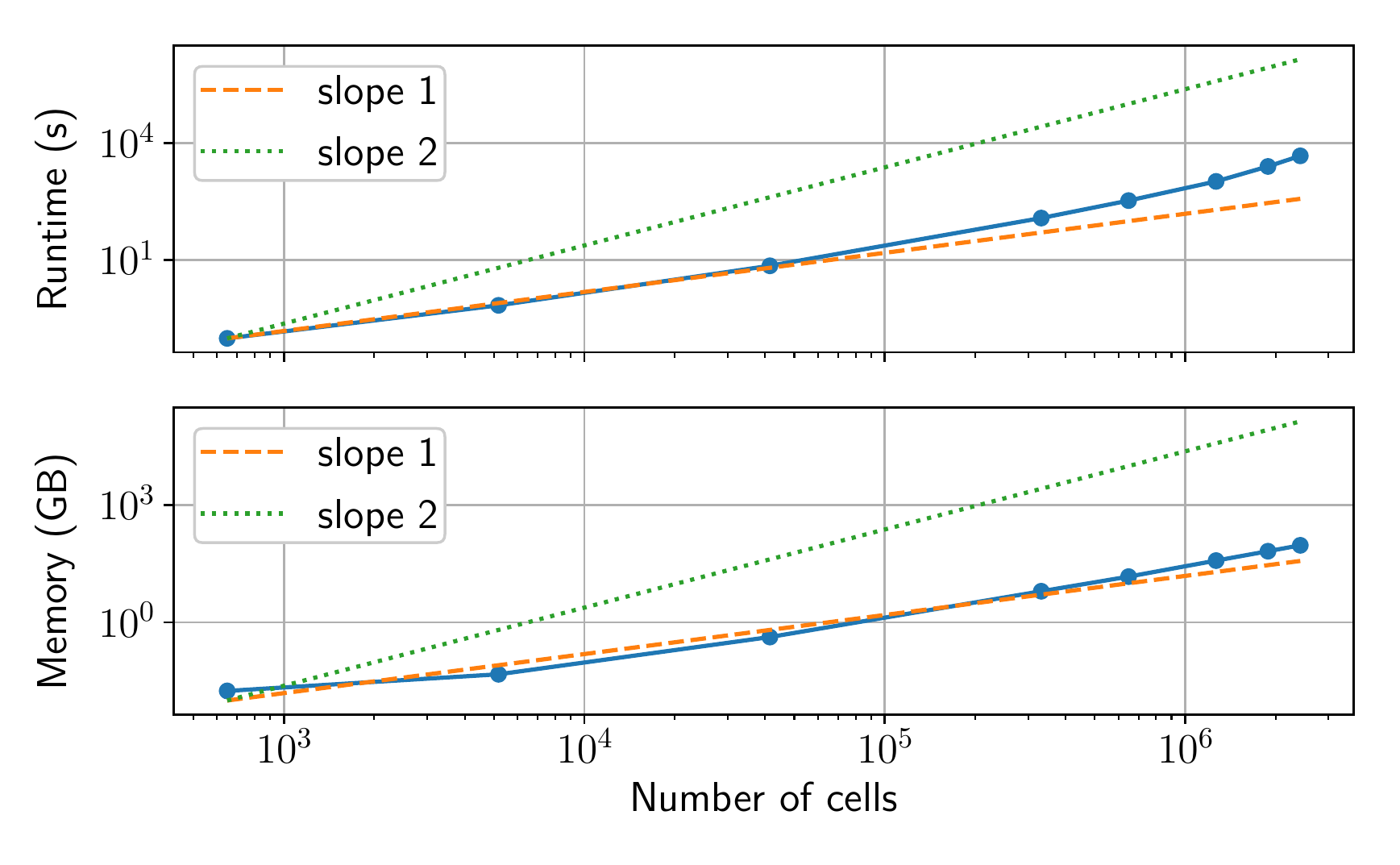}
    \caption{Runtime and memory usage versus mesh size for the Helmholtz Resonator case.}
    \label{fig:reson_mesh_scaling}
\end{figure}

We studied the performance and memory consumption of the Helmholtz resonator case as described in \cref{sec:resonator}. This is a cold flow case without reactions, so we use the \texttt{EPS} module of SLEPc to solve the eigenvalues. We run this case in serial with increasing mesh size. The overall runtime (top plot) and memory consumption (bottom plot) is presented in \cref{fig:reson_mesh_scaling}. As we see, the runtime increases linearly with mesh size until about  $10^5$ cells, and beyond that it slowly starts moving toward a quadratic growth. Most of this computation time in spent in MUMPS \cite{Amestoy:2000:MUMPS, Amestoy:2001:FAM} to compute the matrix factorization, whose cost for sparse matrices typically scales as $\mathcal{O}(N^2)$, so the trend we are seeing here is consistent. The memory consumption for this case shows almost a linear trend throughout. The maximum mesh size that we could run for this case was around 2.4M (which consumed around 96 GB of memory), beyond which we ran out of memory on the system.

\begin{figure}
    \centering
    \includegraphics[width=400pt]{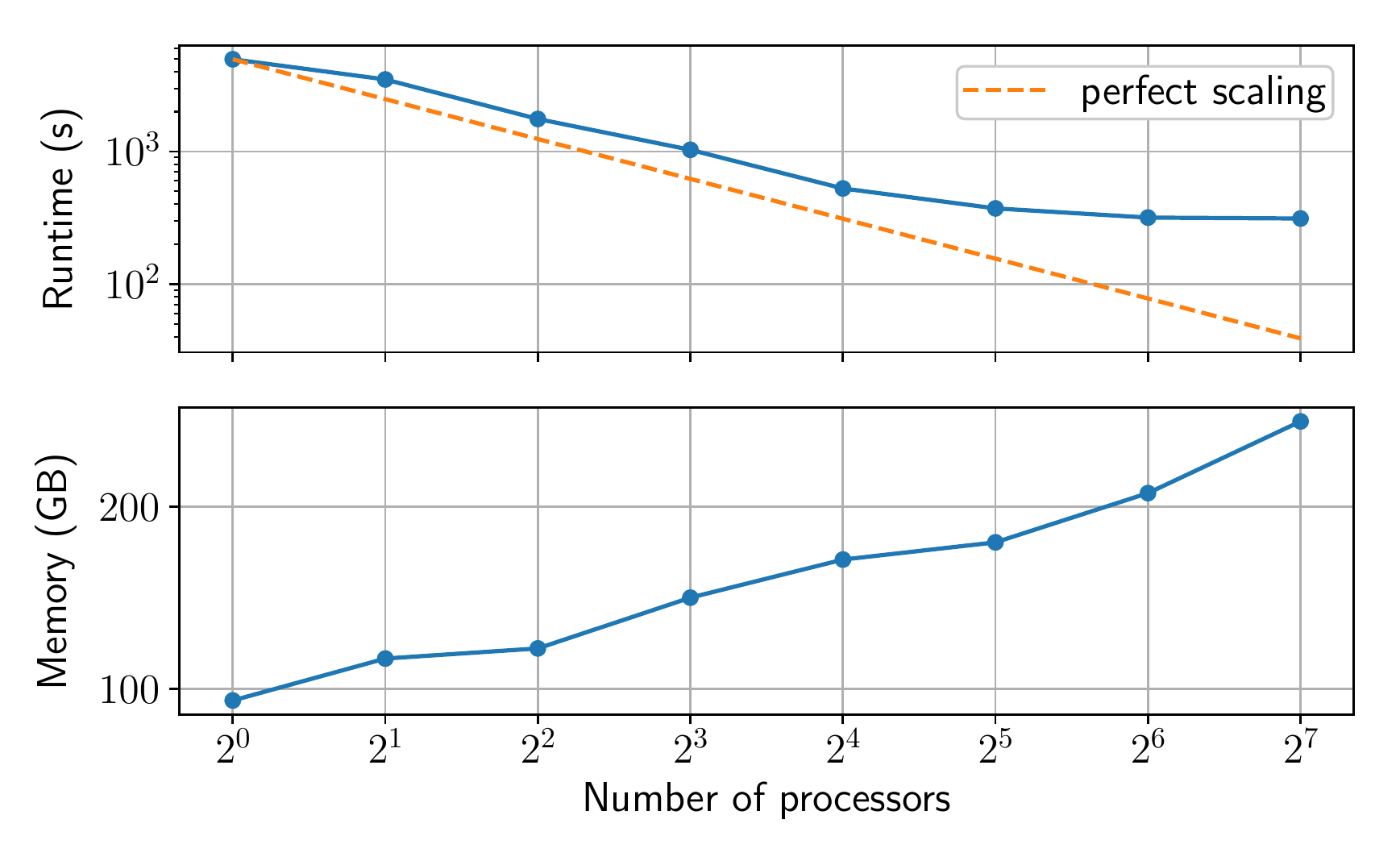}
    \caption{Runtime and memory usage scaling with number of processors for the Helmholtz Resonator case with 2.4M cells.}
    \label{fig:reson_procs_scaling}
\end{figure}

We took the 2.4M mesh case and performed some parallel scaling studies by running the case on increasing number of processors on the aforementioned cluster where each node has 180 GB memory and 40 processors. As needed we spread the processes across multiple nodes to ensure enough memory is always available to run the computation. The parallel scaling results are presented in \cref{fig:reson_procs_scaling}. As we see until about 16 cores we get nearly perfect scaling, beyond which the runtime appears to flatten out as the number of cells per process decrease and the parallel communication costs presumably start dominating. For the fixed mesh size, as we increase the number of processors, the memory demand also increases as we see in the bottom plot in \cref{fig:reson_procs_scaling}. This increase in memory demand is mainly due to matrix data partitioning and duplication of some of the data across all the processors.

\subsection{BH3P burner}
\begin{figure}
    \centering
    \includegraphics[width=400pt]{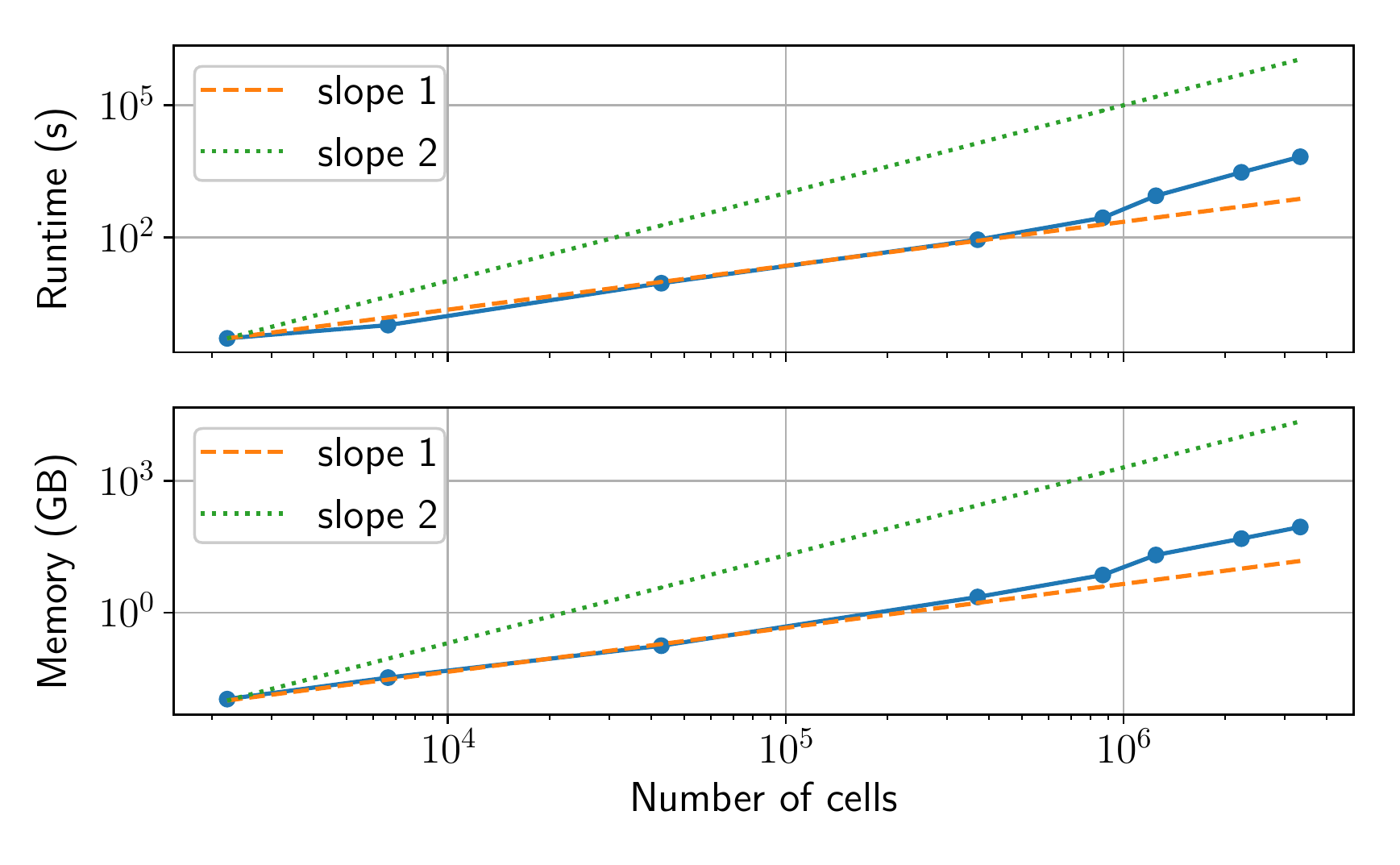}
    \caption{Runtime and memory usage versus mesh size for the Helmholtz Resonator case.}
    \label{fig:bh3p_mesh_scaling}
\end{figure}

\begin{figure}
    \centering
    \includegraphics[width=400pt]{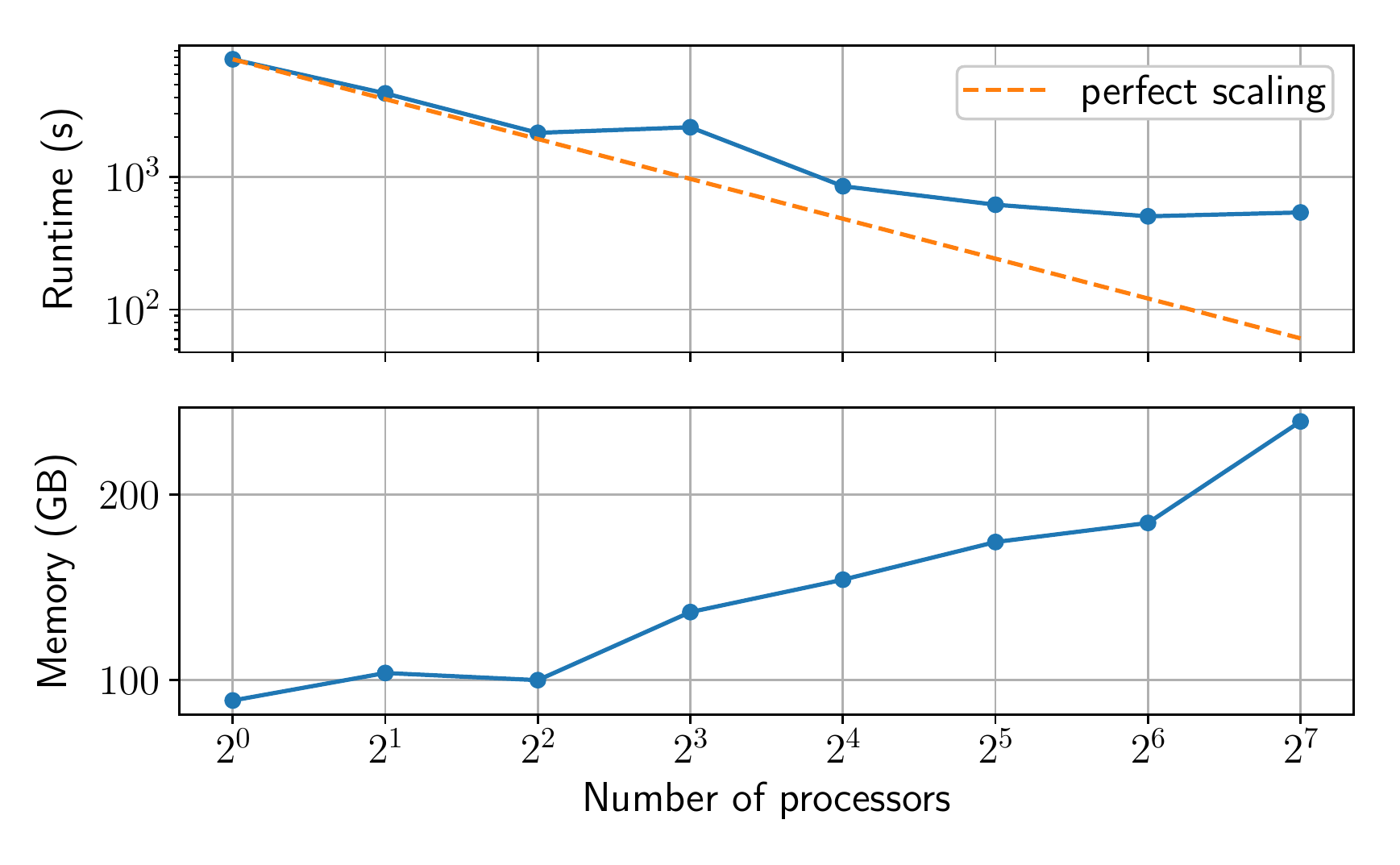}
    \caption{Runtime and memory usage scaling with number of processors for the BH3P burner case with 3.3M cells.}
    \label{fig:bh3p_procs_scaling}
\end{figure}

We do similar scaling studies as in the previous section for the BH3P case as described in \cref{sec:bh3p}. This is a reacting case so we use the \texttt{NEP} module (NLEIGS method) of SLEPc to compute the eigenvalues. We first perform scaling studies with increasing mesh size. The results are presented in \cref{fig:bh3p_mesh_scaling}. As with the Helmholtz resonator case, we first see a linear increase in runtime until about $10^6$ cells, and beyond that a quadratic increase. The memory consumption increases linearly. The parallel scaling results for this case using a 3.3M mesh case are shown in \cref{fig:bh3p_procs_scaling}. In this case again, we see good scaling until about 16 processors, beyond which the runtime starts to flatten out. As seen in the previous case, here again the memory demand increases gradually as we increase the number of processors.

It is encouraging to see that the nonlinear eigenvalue solver (\texttt{NEP}) shows similar trend as the linear eigenvalue solver (\texttt{EPS}), which means that we can solve nonlinear eigenvalue problems as efficiently as linear problems using SLEPc and MUMPS.

\section{Conclusions}\label{sec:conclusions}
In this work we have presented an accurate and efficient way of solving the acoustic Helmholtz eigenvalue problem using the NLEIGS algorithm \cite{Guttel:2014:NCF} implemented in SLEPc's \texttt{NEP} package. We have shown that NLEIGS algorithm works better than some of the earlier proposed linearized iterative algorithms \cite{Nicoud:2007:AIAA}. We have validated the acoustic modal solver implementation in STAR-CCM+ using SLEPc for various cases and have demonstrated its accuracy and efficiency. We have also shown that the SLEPc package along with MUMPS can be used to solve large cases efficiently both in serial and parallel on multiple processors.

In the near future we plan to investigate more thoroughly the SLEPc's CISS method and its different variants and other contour integration methods suggested in \cite{Buschmann:2020:JEGTP, Mohamed:2020:EABM} and compare all these with the NLEIGS method. We also plan to look at more efficient usage of MUMPS \cite{Amestoy:2006:PC} in parallel. We currently invoke MUMPS via SLEPc and use the default settings, however MUMPS has certain options like a Hybrid (OpenMP+MPI) mode that allows more efficient computation in parallel. We plan to investigate this and other internal MUMPS options to allow us to run even bigger cases more efficiently in parallel.

\section*{Acknowledgements}
VH would like to thank Graham Goldin for providing valuable feedback and comments, Doru Caraeni for sharing ideas on the secondary gradient implementation, and Artyom Pogodin for preparing the test cases used in this study and providing useful suggestions.

\end{document}